\documentclass[11pt]{article}
\usepackage{amsbsy}
\usepackage{bm}
\usepackage{ulem}
\usepackage{amssymb,amsmath,amsfonts,amsthm}
\usepackage{color}

\newtheorem {Proposition}{Proposition}[section]
\newtheorem {Lemma}{Lemma}[section]

\newtheorem {corollary}{Corollary}[section]

\newtheorem {Remark}{Remark}[section]
\newtheorem {Assumption}{Assumption}[section]
\usepackage{amssymb,amsmath}
\usepackage{authblk}
\usepackage{epsfig}
\usepackage{enumerate,graphicx}
\usepackage{amsfonts}
\usepackage{graphicx}
\DeclareGraphicsExtensions{.pdf, .jpg, .png , .bmp, .ps}
\def\build#1_#2^#3{\mathrel{\mathop{\kern 0pt#1}\limits_{#2}^{#3}}}
\def\cvl{\build{\ \longrightarrow\ }_{\nti}^{{\mathcal L}}}

\def\cvp{\build{\ \longrightarrow\ }_{\nti}^{{Pr}}}
\def\nti{n \rightarrow \infty}
\DeclareMathOperator{\R}{\mathbb{R}}
\DeclareMathOperator{\X}{\textbf{X}}

\DeclareMathOperator{\Y}{\textbf{Y}}
\DeclareMathOperator{\A}{\textbf{A}}
\DeclareMathOperator{\B}{\textbf{B}}
\DeclareMathOperator{\C}{\textbf{C}}
\DeclareMathOperator{\D}{\textbf{D}}
\DeclareMathOperator{\E}{\textbf{E}}
\DeclareMathOperator{\M}{\textbf{M}}
\newcommand{\bL}{\mathbb{L}}
\def\cvp{\build{\ \longrightarrow\ }_{\nti}^{{\mathbb P}}}
\def\nti{n \rightarrow \infty}

\newcommand{\gX}{\mathbf{X}}

\newcommand{\gY}{\mathbf{Y}}
\def\Var{\mathop{\rm Var}\nolimits}%

\begin{document}
\title{Asymptotic confidence interval for $R^2$ in multiple linear regression}
\author[1]{J. Dedecker}
\affil[1]{Laboratoire MAP5, UMR CNRS  8145, Université Paris Cité, France
}
\author[2]{O. Guedj}
\affil[2]{Laboratoire LaMME, UMR CNRS 8071, Universit\'e d'Evry Val d'Essonne, Université Paris Saclay, France
}
\author[2]{M.L. Taupin}



\maketitle

\begin{abstract}
Following White's approach of robust multiple linear regression \cite{W}, we give asymptotic confidence intervals for the multiple correlation coefficient $R^2$ under minimal moment conditions. We also give the asymptotic joint distribution of the empirical estimators of the individual $R^2$'s. Through different sets of simulations, we show that the procedure is indeed robust (contrary to the procedure involving the near exact distribution of the empirical estimator of $R^2$ is the multivariate Gaussian case) and can be also applied to count linear regression. Several extensions are also discussed, as well as an application to robust screening.
\end{abstract}



\section{Introduction}
Let $(Y, X^{(1)}, \ldots X^{(p)})$ be a random vector with value in ${\mathbb R}^{p+1}$. The variable $Y$ is the response variable, and $(X^{(1)}, \ldots, X^{(p)})$ is the vector of explanatory variables. Assume that all the variables are square integrable. The multiple correlation between $Y$ and $(X^{(1)}, \ldots, X^{(p)})$ can  be expressed using the orthogonal projection of $Y$ onto the ${\mathbb L}^2$-subspace $V$ generated by the variables $(X^{(0)},X^{(1)}, \ldots, X^{(p)})$, where we denote by $X^{(0)}$ the constant variable $X^{(0)}\equiv {\bf 1}$. Let then $\mbox{Proj}_V(Y)$ be the orthogonal projection of $Y$ onto $V$. The response variable $Y$ is uncorrelated to the vector $(X^{(1)}, \ldots, X^{(p)})$ if and only if 
$\mbox{Proj}_V(Y)={\mathbb E}(Y)$. 

To measure the strength of the correlation between $Y$ and  $(X^{(1)}, \ldots, X^{(p)})$, the usual measure is the coefficient $R^2$ defined by (assuming that $\text{Var}(Y)>0$):
$$
R^2= \frac{\text{Var}( \text{Proj}_{V}(Y))}{\text{Var}(Y)} \, .
$$
The coefficient $R^2$ is often referred as "the population variance-accounted-for effect size". As we see, it is defined as the proportion of the variance  of $Y$ that is explained by the best linear predictor in ${\mathbb L}^2$ based on the variables 
${\bf 1}, X^{(1)}, \ldots , X^{(p)}$.

In this paper, we consider the empirical estimator $\widehat{R}^2$ of the coefficient $R^2$ based on $n$ independent copies
$(Y_i, X_{i}^{(1)}, \ldots ,X_{i}^{(p)})_{1 \leq i \leq n}$ of the vector $(Y, X^{(1)}, \ldots , X^{(p)})$. We prove its consistency, asymptotic normality,  and we give an asymptotic confidence interval for $R^2$. Our approach is a robust one: we will make no additional assumptions on the distribution of the vector $(Y, X^{(1)}, \ldots , X^{(p)})$, except moment assumptions necessary for the Central Limit Theorem (CLT) to hold, and to estimate consistently the limiting variance.

Most of the existing literature on the distribution of $\widehat R^2$ deals with the case where the vector $(Y, X^{(1)}, \ldots X^{(p)})$ is normally distributed. In this context, Fisher (1928) \cite{F} obtained the first exact expression of the density function of $\widehat R^2$, involving a Gauss hypergeometric series. Other exact expressions for the distribution and the cumulative distribution   were given by Lee \cite{L1}. An exact expression of the cumulative distribution in terms of a series of Gamma distributions is given in Muirhead (1982, Theorem 5.2.5) \cite{Muir}. The first two principal terms of the latter development were given earlier by Lee \cite{L1}. Some quantile tables were derived by Lee \cite{L2} from this second-order approximation of the cumulative distribution function. Lee's results were then implemented by Kelley \cite{K} in his \texttt{R} MBESS package. 

In the non-Gaussian case, there are relatively few results. Let us cite Muirhead \cite{Muir} Theorem 5.1.6, who gave the asymptotic normality of $\widehat{R}^2$ assuming that $(Y, X^{(1)}, \ldots X^{(p)})$  has an elliptical distribution. 
 The paper closest to ours is that of Ogasawara  \cite{Og}, who gave an Edgeworth expansion for $\sqrt n(\widehat{R}^2-R^2)$, under some conditions on the distributions of the variables, assuming in particular that all the variables $Y, X^{(1)}, \ldots X^{(p)}$ have a moment of order 8. We will compare our result to that of Ogasawara in more details in Remark \ref{Og} of Section \ref{Sec3}. Let us also mention the results of Ogasawara (2008) \cite{Og2} for partial correlation (see Subsection \ref{Sec6.1} for the extension of our method to the the case of partial correlation). 

 As the previous two paragraphs show, the study of the distribution of $\widehat{R}^2$ is important from a historical point of view, the first work on this question being prior to 1930. It is also important from a practical point of view, since $\widehat{R}^2$ is systematically computed in the case of the multivariate linear model, its interpretation being particularly clear in this case. Obtaining a confidence interval of the true $R^2$ from $\widehat{R}^2$ is therefore a perfectly natural question; it gives numerical precision to the quality of the forecasts obtained from the best linear combination of the explanatory variables. We refer to Kelley's article \cite{K2} for further discussion of the practical applications of these confidence intervals. 
 
 As a referee suggested, some of the results of this article may be of interest to practitioners in the Machine Learning community, as $R^2$ (or some related quantities, such as LMG/PLMV indices - see for instance Gr\"omping \cite{Gr}) is often used in this field. Again, we wish to note here that our approach (inspired by White \cite{W}) is a robust one in the sense that we never assume the existence of a model to describe the relationship between the variable $Y$ and the explanatory variables $(X^{(1)}, \ldots,  X^{(p)})$. This robustness property, which is also highlighted in Section \ref{SecRS} devoted to screening, could naturally fit into the context of SAFE Machine learning, as described for example in Guidici and Raffinetti \cite{GR}.  The concept of SAFE machine learning is taken up and developed in detail in the recent article by Giudici \cite{Giu}. Our robust method could then fall under two principles presented in this article: accuracy in the sense that confidence intervals for $R^2$ make it possible to obtain precision on the predictive capacities of the best linear predictor, sustainability in the sense that we work in a model-free context, which does not require the verification of prior hypotheses. Similarly, confidence intervals on the LMG importance measure (see Gr\"omping \cite{Gr}) in a model-free context would be interesting in terms of the explainability principle (even if the notion of variable importance does not always help to explain or interpret the link between the output and the explanatory variables).

The article is organized as follows: in Section \ref{Sec2}, we recall some earlier results due to White \cite{W} necessary for understanding the rest of the article. In Section \ref{Sec3}  we prove the  consistency and asymptotic normality of $\widehat{R}^2$. We also give an estimator of the limiting variance, which together with asymptotic normality provide an asymptotic confidence interval for $R^2$. In Section \ref{Sec4}, we evaluate the performance of our confidence interval on different sets of simulations, and we compare this performance with that of the confidence interval that assumes that the vector $(Y, X^{(1)}, \ldots X^{(p)})$ is Gaussian. In Section \ref{Sec5}, we give the extension of our procedure to the cases of partial correlation, vector-valued random variables, and dependent ($\alpha$-mixing) random variables. In Section \ref{Sec6}, we present additional results and some possible applications of the techniques developed in Section \ref{Sec3}. In Section \ref{Sec4bis}, we give the asymptotic distribution of the joint law of the estimators of the individual $R^2$'s. In Sections \ref{SecSA} and \ref{SecRS} we briefly discuss the connections between our results and sensitivity analysis or screening methods (in the latter case, a small simulation study is carried out).


\section{Previous known results}\label{Sec2}

In this section, we recall White's results \cite{W} concerning the least squares estimator of the coefficients of the orthogonal projection $\text{Proj}_{V}(Y)= \alpha_0 +\alpha_1 X^{(1)}+ \cdots + \alpha_p X^{(p)}$ (consistency, asymptotic normality and estimation of the limiting covariance matrix). This serves two purposes: it recalls the robust procedure described by White \cite{W} in the linear model (without the assumption of normality or homoscedasticity), and it gives an initial idea of the proofs that will enable us to obtain an asymptotic confidence interval for $R^2$ in Section \ref{Sec3}.

\medskip

Let $V_i$ be the ${\mathbb L}^2$-subspace generated by the variables  $(X^{(0)}_i,X^{(1)}_i, \ldots, X^{(p)}_i)$. Let $\varepsilon=Y- \mbox{Proj}_{V}(Y)$ and $ \varepsilon_i= Y_i- \mbox{Proj}_{V_i}(Y_i)$.
 We can then write
$$
Y_i= \mbox{Proj}_{V_i}(Y_i)+ \varepsilon_i.
$$
where the variables $\varepsilon_i$ are such that ${\mathbb E}(X_{i}^{(j)} \varepsilon_i)=0$ for any $j \in \{0, \ldots , p\}$.

Consider the following assumption:
\begin{Assumption}
\label{A_independent}
Assume that  $ {\bf 1}, X^{(1)}, \ldots,  X^{(p)}$ are linearly independent. By linearly independent, we mean ``as a family of random variables in the space ${\mathbb L}^2$", that is: one cannot find ${\bm \gamma}=(\gamma_0, \gamma_1, \ldots , \gamma_p)$ in ${\mathbb R}^{p+1}$ with $\bm \gamma \neq 0$ such that 
$
\gamma_0 + \gamma_1 X^{(1)} + \cdots + \gamma_p X^{(p)} = 0 
$ almost surely.
\end{Assumption}
Note that Assumption \ref{A_independent} is true iff the distribution of $(X^{(1)}, \ldots,  X^{(p)})$ is not supported on an affine subspace of dimension $k<p$.

Under Assumption \ref{A_independent}
one can uniquely write
$$
Y_i= \alpha_0 +\alpha_1 X_{i}^{(1)}+ \cdots + \alpha_p X_{i}^{(p)} +\varepsilon_i  \, .
$$
Let $\bm{\alpha}_{0:p}=(\alpha_0,\alpha_1,\ldots,\alpha_p)^t\in \mathbb{R}^{p+1}$, and $\X$ be the matrix whose $j$th column is 
$$\X^{(j)}=(X_1^{(j)},\ldots,X_n^{(j)})^t.$$
Let $\mathbf{Y}=(Y_1,\ldots,Y_n)^t$ and $\bm{\varepsilon}=(\varepsilon_1,\ldots,\varepsilon_n)^t$. One can write
$$\mathbf{Y}=\mathbf{X}\bm{\alpha}_{0:p}+\bm{\varepsilon}.$$
Let $\mathbf M$ be the matrix such that  $\mathbf M_{j,\ell}={\mathbb E}(X^{(j)} X^{(\ell)})$. By Assumption (\ref{A_independent}) the matrix $\mathbf M$ is invertible and one can easily show that
$$
\bm{\alpha}_{0:p}= \mathbf M^{-1}({\mathbb E}(Y),{\mathbb E}(X^{(1)}Y) \ldots, {\mathbb E}(X^{(p)}Y))^t \, .
$$

\subsection{Least square estimator of ${\bm{\alpha}}_{0:p}$}
Let 
$$
\widehat{\mathbf M}= \frac 1 n \X^t\X \, .
$$
By the strong law of large numbers
\begin{equation}\label{lgn1}
\lim _{n \rightarrow \infty} \widehat{\mathbf M} = \mathbf M \ \text{almost surely.}
\end{equation}
Under Assumption \ref{A_independent}, $\mathbf M$  is invertible and for large enough $n$, $\X^t \X$ is also invertible, and one can define  
$$
\widehat{\bm{\alpha}}_{0:p}=(\widehat \alpha_0, \ldots, \widehat \alpha_p)^t= (\X^t \X)^{-1}\X^t \mathbf{Y}.
$$
Let  ${\mathbb V}_{p+1}$ be the sub-space of ${\mathbb R}^n$ generated by the  $p+1$ columns $\X^{(0)}, \X^{(1)},  \ldots, \X^{(p)}$ of the matrix  $\X$, and $\Pi_{{\mathbb V}_{p+1}}(\mathbf{Y})$ be the orthogonal projection of $\mathbf{Y}$ on  ${\mathbb V}_{p+1}$ (with respect to the euclidean norm).  We have  
$$
 \widehat \alpha_0\X^{(0)}+\widehat \alpha_1 \X^{(1)} + \cdots + \widehat \alpha_p   \X^{(p)}=\Pi_{{\mathbb V}_{p+1}}(\mathbf{Y}),
$$
and
\begin{equation}\label{est}
\widehat{\bm{\alpha}}_{0:p}-\bm{\alpha}_{0:p}=(\widehat \alpha_0, \ldots, \widehat \alpha_p)^t-(\alpha_0, \ldots, \alpha_p)^t = (\X^t \X)^{-1}\X^t\bm{\varepsilon}= \widehat{\mathbf M} ^{-1} \frac 1 n \X^t\bm{\varepsilon}\, .
\end{equation}

\subsection{Consistency and asymptotic normality of $\widehat{\bm{\alpha}}_{0:p}$}\label{sec2.2}

Following White \cite{W}, we obtain the consistency and asymptotic normality of the least square estimators. The proofs of these results being simple and enlightening, we have chosen to recall them because they may make it easier to understand the proofs of Section \ref{Sec3}.

\begin{Proposition}\label{consistEst}
Under Assumption \ref{A_independent}, $\widehat{\bm{\alpha}}_{0:p}$
converges almost surely to $\bm{\alpha}_{0:p}$ as $n\rightarrow \infty$. 
\end{Proposition}

\noindent {\bf Proof of Proposition \ref{consistEst}.} 
From \eqref{lgn1}, we get that
\begin{equation}\label{lgn2}
\lim _{n \rightarrow \infty} \widehat{\mathbf M}^{-1} = \mathbf M ^{-1}\ \text{almost surely.}
\end{equation}
Starting from \eqref{est} and using  \eqref{lgn2}, it suffices to show that
\begin{equation}\label{lgn3}
\lim _{n \rightarrow \infty} \frac 1 n \X^t \bm{\varepsilon}
= 0\ \text{almost surely.}
\end{equation}
By definition of $\varepsilon$, for any  $j \in \{0, \ldots , p\}$, ${\mathbb E}(X^{(j)} \varepsilon)=0$. Applying the strong law of large numbers to the $j$th coordinate of  $ n^{-1} \X^t \bm{\varepsilon}$, we get
\begin{equation}\label{coor}
\lim _{n \rightarrow \infty}\left (\frac 1 n \X^t \bm{\varepsilon}
\right )_j=\lim _{n \rightarrow \infty}\frac 1 n \sum_{i=1}^n X_{i}^{(j)} \varepsilon_i=0\ \text{almost surely. \quad \quad \quad \quad \qed}
\end{equation}

To prove the asymptotic normality, we need moment assumptions.

\begin{Assumption}
    \label{A_moment_1}
    Assume that 
    ${\mathbb E}({X^{(j)}}^2 \varepsilon^2)< \infty$  for all $j=0,\ldots,p$.
\end{Assumption}

\begin{Proposition}\label{cltP}
Under Assumptions \ref{A_independent} and \ref{A_moment_1}, 
let $\mathbf M({\varepsilon})$ be the matrix defined by   
\begin{eqnarray}
    \label{def_M_epsilon}
\mathbf M(\varepsilon)_{j,\ell}={\mathbb E}(X^{(j)} X^{(\ell)} \varepsilon^2) \ \ \text{for any} \ 0\leq j,\ell\leq p.
\end{eqnarray}Then
$$
\sqrt n (\widehat{\bm{\alpha}}_{0:p}-\bm{\alpha}_{0:p})
\cvl {\mathcal N}_{p+1}(0, \mathbf M^{-1}\mathbf M(\varepsilon)\mathbf M^{-1}).
$$
\end{Proposition}

\begin{Remark} Assume that $\varepsilon^2-{\mathbb E}(\varepsilon^2)$ is orthogonal to  the space $S$ generated by the variables $X^{(j)} X^{(\ell)}$ for all  $0\leq j \leq \ell\leq p$ (which is true, for instance,  if $\varepsilon$ is independent of $(X^{(1)}, \ldots, X^{(p)})$). Then the limiting variance matrix in Proposition \ref{cltP} writes
 $$
 \mathbf M^{-1}\mathbf M(\varepsilon)\mathbf M^{-1}= {\mathbb E}(\varepsilon^2) \mathbf M^{-1} \, ,
 $$
 and is the same as in the case where the conditional distribution of $Y$ given $(X^{(1)}, \ldots, X^{(p)})$ is Gaussian with $\mathrm {Var}(Y| X^{(1)}, \ldots, X^{(p)})=\sigma^2$ (homoscedastic case; in that case ${\mathbb E}(\varepsilon^2)=\sigma^2$). This simple observation enabled White \cite{W} to formulate his famous homoscedasticity test, which in fact consists of testing whether $\varepsilon^2-{\mathbb E}(\varepsilon^2)$ is orthogonal to $S$ ($H_0$ hypothesis) or not.
 \end{Remark}

 \medskip

\noindent {\bf Proof of Proposition \ref{cltP}.} Starting from \eqref{est} and using  \eqref{lgn2}, it suffices to prove that
\begin{equation}\label{clt}
 \frac 1 {\sqrt n} \X^t \bm{\varepsilon}
 \cvl {\mathcal N}_{p+1}(0, \mathbf M(\varepsilon)).
 \end{equation}
 Now \eqref{clt} follows from a direct application of the central limit theorem in ${\mathbb R}^{p+1}$, since for any $0\leq j\leq p$,
 ${\mathbb E}(X^{(j)} \varepsilon)=0$ and since Assumption \ref{A_moment_1} holds.

 \qed

\subsection{Estimation of the limiting covariance matrix}\label{sec:est}
Under Assumption \ref{A_moment_1}, by the strong law of large numbers, for any $j,\ell \in\{0, \ldots , p\}$,
\begin{equation}\label{prelim}
\lim_{n \rightarrow \infty} \frac 1 n \sum_{i=1}^n X_{i}^{(j)}X_i^{(\ell)} \varepsilon^2_i = \mathbf M(\varepsilon)_{j,\ell} \quad \text{almost surely}.
 \end{equation}
Of course, the quantity on right hand is not an estimator of $\mathbf M(\varepsilon)_{j,\ell}$ since the $\varepsilon_i$'s are not observed.
White \cite{W} proposed then to replace the $\varepsilon_i$'s by the residuals.
$$
\widehat \varepsilon_i= Y_i-({\bf X} \widehat{\bm{\alpha}}_{0:p})_i
=Y_i-\widehat \alpha_0 X_{i}^{(0)} - \cdots - \widehat \alpha_p X_{i}^{(p)}. 
$$ 
Let then
$$
\widehat {\mathbf M(\varepsilon)}_{j,\ell}= \frac 1 n \sum_{i=1}^n  X_{i}^{(j)}  X_{i}^{(\ell)}  \widehat \varepsilon_i^2  \, .
$$
To prove the consistency of $\widehat {\mathbf M(\varepsilon)}_{j,\ell}$, we shall assume that Assumption \ref{A_moment_1} holds and that:
\begin{Assumption}
\label{A_moment_2}
${\mathbb E}({X^{(j)}}^4)< \infty$ for all $j \in \{0, \ldots , p\}$.
\end{Assumption}

\begin{Proposition}\label{consist}
Under Assumptions \ref{A_moment_1} and \ref{A_moment_2},
for any $j,\ell \in \{0, \ldots , p\}$,
$$
\lim_{n \rightarrow \infty} \widehat {\mathbf M(\varepsilon)}_{j,\ell}=\mathbf M(\varepsilon)_{j,\ell} \quad \text{almost surely}.
$$
Consequently, under Assumptions \ref{A_independent}-\ref{A_moment_2},  
$$
\lim_{n \rightarrow \infty} \widehat{\mathbf M}^{-1} \widehat{\mathbf M(\varepsilon)} \widehat{\mathbf M}^{-1}= {\mathbf M}^{-1} {\mathbf M(\varepsilon)} {\mathbf M}^{-1} \quad \text{almost surely}.
$$
\end{Proposition}
 
\begin{Remark} The Assumptions \ref{A_moment_1} and \ref{A_moment_2} are equivalent to: ${\mathbb E}({X^{(j)}}^2 Y^2)< \infty$ and ${\mathbb E}({X^{(j)}}^4)< \infty$ for any $j \in \{0, \ldots , p\}$.
\end{Remark}

\begin{Remark}\label{test}
As usual, Propositions \ref{cltP} and \ref{consist} enable to obtain confidence regions or test procedures for the coefficients $\alpha_i$. For instance, if one wants to test
$$
H_0 : \alpha_{i_1}= \cdots = \alpha_{i_k}=0 \quad \text{for some $0\leq i_1 < i_2 < \cdots < i_k \leq p$}
$$
one can proceed as follows. Let $\mathbf C_{k}$ be the matrix with $k$ rows and $p$ columns, such that all coordinates of the $j$th row are 0 except the coordinate $i_j$ which is equal to 1. Then, by Proposition \ref{cltP}, 
$$
\sqrt n ((\widehat \alpha_{i_1}, \ldots , \widehat \alpha_{i_k})^t-(\alpha_{i_1}, \ldots , \alpha_{i_k})^t)
\cvl {\mathcal N}_{k}(0, \mathbf C_k\mathbf M^{-1}\mathbf M(\varepsilon)\mathbf M^{-1} \mathbf C_k^t).
$$
Now, under Assumptions \ref{A_independent}-\ref{A_moment_2}, by applying Proposition \ref{consist} 
$$
\widehat{\mathbf V}_k = \mathbf C_k \widehat{\mathbf M}^{-1}\widehat{\mathbf M(\varepsilon)}\widehat{\mathbf M}^{-1} \mathbf C_k^t \quad \text{converges almost surely to} \quad 
\mathbf C_k\mathbf M^{-1}\mathbf M(\varepsilon)\mathbf M^{-1} \mathbf C_k^t \, .
$$
If $\mathbf M(\varepsilon)$ is invertible, 
then the matrix $\mathbf C_k\mathbf M^{-1}\mathbf M(\varepsilon)\mathbf M^{-1} \mathbf C_k^t$ is also invertible
and  
 $\widehat{\mathbf V}_k^{-1/2}$ is well defined (for $n$ large enough). Let then $(\xi_1, \ldots ,\xi_k)^t=\sqrt n \widehat{\mathbf V}_k^{-1/2}(\widehat \alpha_{i_1}, \ldots , \widehat\alpha_{i_k})^t$. Under $H_0$, we see that the test statistic $\sum_{i=1}^k \xi_i^2$ is such that
 $$
 \sum_{i=1}^k \xi_i^2 \cvl \chi^2(k) \, .
 $$
\end{Remark}

\noindent {\bf Proof of Proposition \ref{consist}.}
Write
$
\widehat \varepsilon^2_i = \varepsilon^2_i +R_i
$, with
$
R_i = 2\varepsilon_i(\widehat \varepsilon_i- \varepsilon_i)+(\widehat \varepsilon_i- \varepsilon_i)^2
$. The following upper bound holds
\begin{equation}\label{RiB1}
|R_i| \leq 2 \sum_{\ell=0}^p |(\widehat \alpha_\ell-\alpha_\ell) X_{i}^{(\ell)} \varepsilon_i| + (p+1) \sum_{\ell=0}^p (\widehat \alpha_\ell-\alpha_\ell)^2 {X_{i}^{(\ell)}}^2 \, .
\end{equation}
Tacking into account \eqref{prelim}, to prove Proposition \ref{consist}, it suffices to show that, for $0\leq j,k\leq p$ and $1\leq i\leq n$,
\begin{equation}\label{inter}
\lim_{n \rightarrow \infty} \frac 1 n \sum_{i=1}^n |X_{i}^{(j)}X_{i}^{(k)} R_i| =0  \quad \text{almost surely}.
\end{equation}
By \eqref{RiB1}, we have 
$$
\frac 1 n \sum_{i=1}^n |X_{i}^{(j)}X_{i}^{(k)} R_i| \leq 
2 \sum_{\ell=0}^p |(\widehat \alpha_\ell-\alpha_\ell)| \frac 1 n \sum_{i=1}^n |X_{i}^{(j)}X_{i}^{(k)} X_{i}^{(\ell)} \varepsilon_i | + (p+1) \sum_{\ell=0}^p (\widehat \alpha_\ell-\alpha_\ell)^2 \frac 1 n \sum_{i=1}^n |X_{i}^{(j)}X_{i}^{(k)} {X_{i}^{(\ell)}}^2 | \, .
$$
Since $|\widehat \alpha_\ell-\alpha_\ell|$ converges almost surely to 0, \eqref{inter} will be satisfied provided, for $0\leq j,k,\ell\leq p,$
$$
\limsup_{n \rightarrow \infty} \frac 1 n \sum_{i=1}^n |X_{i}^{(j)}X_{i}^{(k)} X_{i}^{(\ell)} \varepsilon_i | < \infty \quad \text{and}  \quad \limsup_{n \rightarrow \infty} \frac 1 n \sum_{i=1}^n |X_{i}^{(j)}X_{i}^{(k)} {X_{i}^{(\ell)}}^2 | < \infty \quad \text{almost surely},
$$
which is true by applying once again the strong law of large numbers, since  under Assumptions \ref{A_moment_1} and \ref{A_moment_2} we have 
${\mathbb E}(|X^{(j)}X^{(k)} X^{(\ell)} \varepsilon |)< \infty$ and ${\mathbb E}(|X^{(j)}X^{(k)} {X^{(\ell)}}^2 |)< \infty$ for any $0\leq j,k,\ell\leq p$.
\qed

\section{Asymptotic distribution and confidence interval for $R^2$}\label{Sec3}
\setcounter{equation}{0}

In this Section   we prove the  consistency and asymptotic normality of $\widehat{R}^2$ (see Proposition \ref{limitlawR2}). We also give an estimator of the limiting variance, which together with asymptotic normality provide an asymptotic confidence interval for $R^2$ (see Corollary \eqref{CI}).

\medskip

Assuming that $\text{Var}(Y)>0$ (which means that the random variable $Y$ is not almost surely constant), the multiple correlation coefficient $R^2$ is defined by:
$$
R^2= \frac{\text{Var}( \text{Proj}_{V}(Y))}{\text{Var}(Y)} \, .
$$
For  $\varepsilon_i= Y_i - \text{Proj}_{V_i}(Y_i)$, recall that under Assumption \ref{A_independent} one can uniquely write
$$
Y_i= \text{Proj}_{V_i}(Y_i)+\varepsilon_i=\alpha_0 {\mathbf 1}+ \alpha_1 X_{i}^{(1)} + \cdots + \alpha_p X_{i}^{(p)} +\varepsilon_i .
$$

For any vector $\mathbf z=(z_1, \ldots, z_n)^t$ in ${\mathbb R}^n$, the  empirical variance of $\mathbf z$
is given by $$\text{Var}_n(\mathbf z)= \frac 1 n \sum_{i=1}^n (z_i - \overline z)^2.$$ 
Let ${\bf 1}_n= (1, \ldots, 1)^t$ and let $\Pi_{E}$ be the orthogonal projection (with respect to the Euclidean norm) on a subspace $E$ of ${\mathbb R}^n$. Let ${\mathbb V}_{p+1}$ be the subspace of
${\mathbb R}^n$ generated by $({\bf 1}_n, \X^{(1)}, \dots , \X^{(p)})$.
The empirical estimator $\widehat R^2$  of $R^2$ is then
\begin{equation}\label{defR2C}
\widehat R^2= \frac{\text{Var}_n( \Pi_{{\mathbb V}_{p+1}}(\Y))}{\text{Var}_n(\Y)} \, .
\end{equation}
Our first goal is to give a simple expression of this estimator. Let $\bm \alpha_{1:p}= (\alpha_1, \ldots, \alpha_p)^t$ and $\widehat {\bm \alpha}_{1:p}= (\widehat\alpha_1, \ldots, \widehat \alpha_p)^t$; for the sake of simplicity, we shall omit the indexes and write $\bm \alpha = \bm \alpha_{1:p}$ and $\widehat{\bm \alpha} = \widehat {\bm \alpha}_{1:p}$. Let also $\bm \theta = (\theta_1, \ldots , \theta_p)^t $ and 
$\widehat{\bm{\theta}}=(\widehat \theta_1, \ldots, \widehat \theta_p)^t$, where
\begin{equation}\label{def:theta}
\theta_k= \frac{\text{Cov}(Y, X^{(k)})}{\text{Var}(Y)} 
\quad \text{and} \quad 
\widehat \theta_k =\frac{\sum_{i=1}^n (Y_i-\overline \Y)(X_{i}^{(k)}-\overline {\X}^{(k)})}{\sum_{i=1}^n (Y_i-\overline \Y)^2}\, .
\end{equation}
The following lemma  gives the expression of $R^2$ and $\widehat R^2$
in terms of ${\bm{\theta}}$, ${\bm{\alpha}}$, ${\widehat{\bm{\theta}}}$ and 
${\widehat{\bm{\alpha}}}$.

\begin{Lemma}\label{scalarproduct}
Under Assumption \ref{A_independent},
$
 R^2= {\bm{\theta}}^t{\bm{\alpha}}
$
and
$
\widehat R^2= {\widehat{\bm{\theta}}}^t{\widehat{\bm{\alpha}}}
$.
\end{Lemma}

\noindent{\bf Proof of Lemma \ref{scalarproduct}.} We prove the second point, the first point being proved in the same way. 
Write
$$
\Pi_{{\mathbb V}_{p+1}}(\Y)- \overline{\Pi_{{\mathbb V}_{p+1}}}(
\Y){\bf 1}_n= \widehat \alpha_1 (\X^{(1)}-\overline{ \X}^{(1)}{\bf 1}_n) + \cdots 
+ \widehat \alpha_p (\X^{(p)}-\overline \X^{(p)}{\bf 1}_n) \, .
$$
Now, if $S_p$ is the subspace of 
${\mathbb R}^n$ generated by $((\X^{(1)}-\overline \X^{(1)}{\bf 1}_n)), \dots , (\X^{(p)}-\overline \X^{(p)}{\bf 1}_n)))$ it is easy to check that $$\widehat \alpha_1 (\X^{(1)}-\overline \X^{(1)}{\bf 1}_n)) + \cdots 
+ \widehat \alpha_p (\X^{(p)}-\overline \X^{(p)}{\bf 1}_n))= \Pi_{S_p}(\Y-\overline{\Y}{\bf 1}_n)
.$$ 

Since  by Assumption \ref{A_independent} the variables  ${\bf 1}, X^{(1)}, \ldots, X^{(p)}$ are linearly independent in ${\mathbb L}^2$, the variance matrix $\mathbf M_0$ of $(X^{(1)}, \ldots, X^{(p)})^t$ defined by 
\begin{equation}\label{defM0}
(\mathbf M_0)_{j,\ell}= \text{Cov}(X^{(j)},X^{(\ell)}) \quad  j,\ell \in \{1, \ldots, p\}
\end{equation}
is invertible. Denote by ${\mathbb X}_0$ the matrix defined 
by 
\begin{equation}\label{defX0}
{\mathbb X}_0=((X_{i}^{(j)}-\overline \X^{(j)}))_{1 \leq i \leq n, 1 \leq j \leq p} \, . 
\end{equation}
By the strong law of large numbers 
$$
\lim_{n \rightarrow \infty} \frac 1 n {\mathbb X}_0^t {\mathbb X}_0=\mathbf M_0 \quad \text{almost surely}.
$$
Since $\mathbf M_0$ defined in (\ref{defM0}) is invertible, it follows that ${\mathbb X}_0^t {\mathbb X}_0$ is also invertible for $n$ large enough, in such a way that 
\begin{equation}\label{def:alpha}
\widehat{\bm{\alpha}}=
({\mathbb X}_0^t {\mathbb X}_0)^{-1}{\mathbb X}_0^t  (\mathbf{Y}-\overline{\mathbf{Y}}{\bf 1}_n),
\end{equation}
and
$$ \Pi_{S_p}((\Y-\overline \Y{\bf 1}_n))
={\mathbb X}_0({\mathbb X}_0^t {\mathbb X}_0)^{-1}{\mathbb X}_0^t 
(\Y-\overline \Y{\bf 1}_n).
$$
Let $\|\cdot \|_{\text{eu}}$ be the Euclidean norm on ${\mathbb R}^n$. The following equality hold:
\begin{align}
n\text{Var}_n( \Pi_{{\mathbb V}_{p+1}}(
\Y))
&=\left \| \Pi_{S_p}(\Y-\overline \Y{\bf 1}_n)
\right \|_{\text{eu}}^2
\nonumber \\
&= 
(\Y-\overline \Y{\bf 1}_n)^t
{\mathbb X}_0({\mathbb X}_0^t {\mathbb X}_0)^{-1}{\mathbb X}_0^t {\mathbb X}_0
({\mathbb X}_0^t {\mathbb X}_0)^{-1}{\mathbb X}_0^t (\Y-\overline \Y{\bf 1}_n)
\nonumber  \\
&=(\Y-\overline \Y{\bf 1}_n)^t{\mathbb X}_0({\mathbb X}_0^t {\mathbb X}_0)^{-1}
{\mathbb X}_0^t (\Y-\overline \Y{\bf 1}_n)\nonumber  \\
&= (\Y-\overline \Y{\bf 1}_n)^t{\mathbb X}_0 
\widehat{\bm{\alpha}}
\, .\label{step1}
\end{align}
The result follows from \eqref{defR2C} and  \eqref{step1} 
since
\begin{equation}\label{step2}
  \frac{(\Y-\overline \Y{\bf 1}_n)^t{\mathbb X}_0}{n\text{Var}_n(\Y)}=
  \widehat{\bm{\theta}}^t
  \,  .   
\end{equation}
\hfill $\square$

\medskip

Now, to build a confidence interval for $R^2$, we need to describe the asymptotic distribution of 
$$
\sqrt{n}\Big([ \widehat{\bm{\alpha}}:\widehat{\bm{\theta}}]
-[ {\bm{\alpha}}:{\bm{\theta}}]\Big):=
  \sqrt n \left ( (\widehat \alpha_1, \ldots, \widehat \alpha_p, \widehat \theta_1, \ldots, \widehat \theta_p)^t-(\alpha_1, \ldots, \alpha_p, \theta_1, \ldots,  \theta_p)^t\right ) \, .
$$
Let us then define the two matrices involved in this asymptotic distribution. For $k \in \{1, \ldots, p\}$, let 
\begin{equation}\label{def:ek}
e^{(k)} = X^{(k)}-\mbox{Proj}_{W}(X^{(k)}) \quad \text{and} \quad e_{i}^{(k)} = X_{i}^{(k)}-\mbox{Proj}_{W_i}(X_{i}^{(k)}),
\end{equation}
where $W$ is the subspace of ${\mathbb L}^2$ generated by ${\bf 1}$ and $Y$, and $W_i$ is the subspace of ${\mathbb L}^2$ generated by ${\bf 1}$ and $Y_i$. Let also ${\bm e}^{(k)}=(e_{1}^{(k)}, \ldots, e_{n}^{(k)})^t $.

\begin{Assumption}
    \label{A_moment_3}
    ${\mathbb E}((X^{(j)}-{\mathbb E}(X^{(j)}))^2 \varepsilon^2)< \infty$ and 
    ${\mathbb E}((Y-{\mathbb E}(Y))^2 {e^{(j)}}^2)< \infty \mbox{  for all }j \in \{1, \ldots , p\}.$
\end{Assumption}
 Let then $\A$ be the $2p\times 2p$ symmetric matrix defined as follows: 
\begin{itemize}
\item If $(j,k) \in \{1, \ldots, p\}^2$ then 
$$\A_{j,k}={\mathbb E}((X^{(j)}-{\mathbb E}(X^{(j)}))(X^{(k)}-{\mathbb E}(X^{(k)}))\varepsilon^2);$$ 
\item If $(j,k) \in \{p+1, \ldots, 2p\}^2$ then $$\A_{j,k}={\mathbb E}((Y-{\mathbb E}(Y))^2 e^{(j-p)} e^{(k-p)});$$ \item  $(j,k) \in \{1, \ldots , p\} \times \{p+1, \ldots, 2p\}$ then $$\A_{j,k}={\mathbb E}((X^{(j)}-{\mathbb E}(X^{(j)}))(Y-{\mathbb E}(Y)) \varepsilon e^{(k-p)}).$$
\end{itemize}
Let $\delta_{j,k}=0$ if $j\not=k$ and $\delta_{j,j}=1$ and let $\B$ be the $2p\times 2p$ symetric matrix defined as follows: 
\begin{itemize}
    \item 
If $(j,k) \in \{1, \ldots, p\}^2$ then 
$\B_{j,k}=(\mathbf M_0^{-1})_{j,k}$ (see  \eqref{defM0} for the definition of $\mathbf M_0$);
\item If $(j,k) \in \{p+1, \ldots, 2p\}^2$ then $\B_{j,k}=(\text{Var}(Y))^{-1}\delta_{j,k}$; 
\item If $(j,k) \in \{1, \ldots , p\} \times \{p+1, \ldots, 2p\}$ then $\B_{j,k}=0$. 
\end{itemize}

\begin{Proposition}\label{limitlaw}
Under Assumptions \ref{A_independent} and \ref{A_moment_3} 
$$
  \sqrt{n}\Big([ \widehat{\bm{\alpha}}:\widehat{\bm{\theta}}]-
  [ {\bm{\alpha}}:{\bm{\theta}}]\Big)
  \quad \cvl
  {\mathcal N}_{2p}(0,{\mathbf B} {\mathbf A} {\mathbf B}).
$$
\end{Proposition}
\noindent{\bf Proof of Proposition \ref{limitlaw}.} 
Let $\widehat \B$ be the $2p \times 2p$ symetric matrix defined as follows: 
\begin{itemize}
    \item If $(j,k) \in \{1, \ldots, p\}^2$ then 
$\widehat \B_{j,k}=  n ({\mathbb X}_0^t {\mathbb X}_0)^{-1}_{j,k}$;
\item If $(j,k) \in \{p+1, \ldots, 2p\}^2$ then $\widehat \B_{j,k}=(\text{Var}_n(\Y))^{-1}\delta_{j,k}$; 
\item If $(j,k) \in \{1, \ldots , p\} \times \{p+1, \ldots, 2p\}$ then $\widehat \B_{j,k}=0$. 
\end{itemize}
Starting from \eqref{def:theta} and \eqref{def:alpha}, and noting that 
\begin{align*}
  Y_i - \overline{\mathbf{Y}}&= \alpha_1(X_{i}^{(1)}-\overline{\X}^{(1)}) + \cdots + \alpha_p (X_{i}^{(p)}- \overline{\X}^{(p)})+ (\varepsilon_i- \overline{\bm{\varepsilon}})\, , \\
  (X_{i}^{(k)}- \overline {\mathbf{X}}^{(k)})&=\theta_k (Y_i - \overline{\mathbf{Y}}) +(e_{i}^{(k)}- \overline{\bm e}^{(k)})\, \quad \text{for $k \in \{1, \ldots, p \}$\, ,}
\end{align*}
we see that
\begin{equation}\label{step1bis}
\sqrt{n}\left([ \widehat{\bm{\alpha}}:\widehat{\bm{\theta}}]-
  [ {\bm{\alpha}}:{\bm{\theta}}]\right)= \widehat \B \frac{1}{\sqrt n} (\bm{\varepsilon}^t
{\mathbb X}_0, (\Y-\overline \Y{\bf 1}_n)^t \mathbf E)^t \, ,
\end{equation}
where $\E$ is the $n\times p$ matrix such that, for $(i,j) \in \{1, \ldots, n\} \times \{1, \ldots p\}$, 
$
(\E_{i,j})=e_{i}^{(j)}.
$ 

Denote by $\widetilde{\mathbb X}_0$ the matrix defined 
by 
\begin{equation}\label{defX0tilde}
\tilde{\mathbb X}_0=(X_{i}^{(j)}-{\mathbb E}(X^{(j)}))_{1 \leq i \leq n, 1 \leq j \leq p} \, . 
\end{equation}
One can easily check that
\begin{equation}\label{step2bis_ext}
\frac{1}{\sqrt n} \left ((\bm{\varepsilon}^t
{\mathbb X}_0, (Y_1- \overline{\mathbf{Y}}, \ldots, Y_n-\overline{\mathbf{Y}}) \mathbf E)^t -
(\bm{\varepsilon}^t
\tilde {\mathbb X}_0, (Y_1- {\mathbb E}(Y), \ldots, Y_n-{\mathbb E}(Y)) \mathbf E)^t\right ) \cvp 0\, .
\end{equation}
Now, by the multivariate central limit theorem,
\begin{equation}\label{step3}
\frac{1}{\sqrt n} 
(\bm{\varepsilon}^t
\tilde {\mathbb X}_0, (Y_1- {\mathbb E}(Y), \ldots, Y_n-{\mathbb E}(Y)) \mathbf E)^t \cvl \mathcal N_{2p}(0, \A).
\end{equation}
The result follows from \eqref{step1bis}, \eqref{step2bis_ext} and \eqref{step3}, since  by the strong law of large numbers,
\begin{equation}\label{step4}
    \lim_{n\rightarrow \infty} \widehat {\mathbf B} = {\mathbf B} \quad \text{almost surely.}
\end{equation}
 \hfill $\square$

\medskip

As a consequence, we get 
\begin{Proposition}\label{limitlawR2}
Under Assumptions \ref{A_independent} and \ref{A_moment_3} 
$$
  \sqrt n \left ( \widehat R^2 - R^2\right )  \cvl {\mathcal N}(0, V).
$$
where 
\begin{equation}\label{defV}
V= (  \theta_1, \ldots,  \theta_p,  \alpha_1, \ldots,  \alpha_p){\mathbf B}{\mathbf A}{\mathbf B}(  \theta_1, \ldots,  \theta_p,  \alpha_1, \ldots, \alpha_p)^t \, .
\end{equation}
\end{Proposition}
\begin{Remark}\label{Og}
    Ogasawara \cite{Og} has given an Edgeworth expansion of $\sqrt n  ( \widehat R^2 - R^2 ) $, which is a more precise result than Proposition \ref{limitlawR2}, but requires some conditions on the distribution of $(Y, X^{(1)}, \cdots, X^{(p)})$. In particular, he required that the all the variables have a moment of order 8. Note that, if we do not look for an Edgeworth expansion, the method of Ogasawara consists in expressing $R^2$ as a differentiable function of $S=((\mathrm{Cov}(X^{(i)},X^{(j)}))_{1\leq i \leq j \leq p},(\mathrm{Cov}(Y,X^{(i)}))_{1\leq i  \leq p}, \mathrm{Var}(Y))$, proving the CLT for $\sqrt{n}(\widehat S -S)$ ($\widehat S$ being the empirical estimator of $S$), and applying the delta method. The differences with our approach are the following :
    \begin{itemize}
        \item To prove the CLT for $\sqrt{n}(\widehat S -S)$, one needs moments of order 4 for all variables, which is a more restrictive condition than the moment conditions of our Proposition \ref{limitlawR2}.
        \item Applying the delta method to a function of $S$ implies that the limiting variance $V$ will be expressed as a function of a $q \times q$ matrix, where $q=(p+1)(p+2)/2$ (the matrix ${\bf \Omega}$ in \cite{Og}), while we can express $V$ as a function of a $2p \times 2p$ matrix (the matrix ${\mathbf B}{\mathbf A}{\mathbf B}$, see \eqref{defV}). Note that, with our expression of $V$ we are able to give simple sufficient conditions ensuring  that $V>0$ (see Lemma \ref{V} below). 
        \item Finally our approach can be easily extended to the the case where the variables $Y, X^{(1)}, \cdots, X^{(p)}$ are vector-valued (see Section \ref{sec:ext}).
    \end{itemize}
\end{Remark}
\noindent{\bf Proof of Proposition \ref{limitlawR2}.} Starting from Lemma \ref{scalarproduct} and Proposition \ref{limitlaw}, it suffices to apply the delta-method to the function $\phi: {\mathbb R}^{p} \times {\mathbb R}^{p}\rightarrow {\mathbb R} $ defined by 
$$
\phi(x,y)=x^ty \, .
$$
The  proof Proposition \ref{limitlawR2} is complete by evaluating the differential $D\phi_{x,y}$ of $\phi$ at point $(x,y)$: 
$$
D\phi_{x,y}(h_1, h_2)=y^th_1+x^th_2 = (y,x)^t(h_1,h_2) \, . 
\quad \quad \quad \quad \quad \square
$$

\medskip

To build a confidence interval for $R^2$, it remains to find a consistent estimator of $V$. We shall simply replace each element in the definition of $V$ by its empirical counterpart.

Let then $\widehat {\mathbf A}$ be the $2p \times 2p$ symmetric matrix defined as follows: 
\begin{itemize}
\item If $(j,k) \in \{1, \ldots, p\}^2$ then 
$$\widehat {\mathbf A}_{j,k}=\frac 1 n \sum_{i=1}^n (X_{i}^{(j)}-\overline{\mathbf{X}}^{(j)})(X_{i}^{(k)}-\overline{\mathbf{X}}^{(k)})\widehat \varepsilon^2_i;$$ 
\item If $(j,k) \in \{p+1, \ldots, 2p\}^2$ then $$\widehat {\mathbf A}_{j,k}= \frac 1 n \sum_{i=1}^n (Y_i-\overline{\mathbf{Y}})^2 \widehat e^{(j-p)}_{i} \widehat e^{(k-p)}_{i},$$ 
where $\widehat e^{(k-p)}_{i}=(X^{(k-p)}_i- \overline {\mathbf{X}}^{(k-p)}) - \widehat \theta_{k-p} (Y_i -\overline{\mathbf{Y}})$;
\item If $(j,k) \in \{1, \ldots , p\} \times \{p+1, \ldots, 2p\}$ then $$\widehat {\mathbf A}_{j,k}= \frac 1 n \sum_{i=1}^n(X_{i}^{(j)}- \overline {\mathbf{X}}^{(j)})(Y_i-\overline{\mathbf{Y}}) \widehat \varepsilon_i \widehat e^{(k-p)}_{i}.$$
\end{itemize}

\begin{Proposition}\label{estim}
 Let 
\begin{equation}\label{hatV}
\widehat V_n= ( \widehat  \theta_1, \ldots,  \widehat \theta_p,  \widehat \alpha_1, \ldots,  \widehat \alpha_p)\widehat {\mathbf B} \widehat {\mathbf A}\widehat {\mathbf B}(  \widehat \theta_1, \ldots, \widehat  \theta_p,  \widehat  \alpha_1, \ldots,  \widehat \alpha_p)^t \, .
\end{equation}
Under Assumptions \ref{A_independent} and \ref{A_moment_2} and  ${\mathbb E}(Y^4)< \infty$, we have
$$
\lim_{\nti} \widehat V_n =V \quad \mbox{almost surely.} 
$$
\end{Proposition}

\noindent{\bf Proof of Proposition \ref{estim}.} Recall that, by the strong law of large numbers, 
$$
   \lim_{n \rightarrow \infty}(\widehat \alpha_1, \ldots, \widehat \alpha_p, \widehat \theta_1, \ldots, \widehat \theta_p)^t=(\alpha_1, \ldots, \alpha_p, \theta_1, \ldots,  \theta_p)^t \quad \text{almost surely,}
$$
and $$\lim_{\nti}\widehat {\mathbf B} ={\mathbf B} \quad \text{ almost surely.}$$
Hence, it remains to prove that $\widehat {\mathbf A}$ converges to ${\mathbf A}$ almost surely.

To prove this point, we first introduce the matrix $\widetilde {\mathbf A}$ defined as follows:
\begin{itemize}
\item 
  If $(j,k) \in \{1, \ldots, p\}^2$ then 
$$\widetilde {\mathbf A}_{j,k}=\frac 1 n \sum_{i=1}^n (X_{i}^{(j)}-\overline{\mathbf{X}}^{(j)})(X_{i}^{(k)}-\overline{\mathbf{X}}^{(k)}) \varepsilon^2_i;$$ 
\item If $(j,k) \in \{p+1, \ldots, 2p\}^2$ then $$\widetilde {\mathbf A}_{j,k}= \frac 1 n \sum_{i=1}^n (Y_i-\overline{\mathbf{Y}})^2  e^{(j-p)}_{i} e^{(k-p)}_{i};$$ 
\item If $(j,k) \in \{1, \ldots , p\} \times \{p+1, \ldots, 2p\}$ then $$\widetilde {\mathbf A}_{j,k}= \frac 1 n \sum_{i=1}^n(X_{i}^{(j)}- \overline {\mathbf{X}}^{(j)})(Y_i-\overline{\mathbf{Y}})  \varepsilon_i   e^{(k-p)}_{i}.$$
\end{itemize}
Since ${\mathbb E}(Y^4)< \infty$ and   ${\mathbb E}({X^{(j)}}^4)< \infty$, we  deduce from the strong law of large numbers that $$\lim_{\nti}\widetilde {\mathbf A}={\mathbf A} \quad \text{almost surely.}$$

Let us now prove that, for  $(j,k) \in \{1, \ldots, p\}^2$, $$
\lim_{\nti}\vert \widehat {\mathbf A}_{j,k}- \widetilde {\mathbf A}_{j,k}\vert = 0 \quad \text{almost surely.}$$
We write
$
\widehat \varepsilon^2_i = \varepsilon^2_i +R_i
$, with
$
R_i = 2\varepsilon_i(\widehat \varepsilon_i- \varepsilon_i)+(\widehat \varepsilon_i- \varepsilon_i)^2
$. One can easily see that 
\begin{equation}\label{B1}
|R_i| \leq 2 \sum_{\ell=0}^p |(\widehat \alpha_\ell-\alpha_\ell) X_{i}^{(\ell)} \varepsilon_i| + (p+1) \sum_{\ell=0}^p (\widehat \alpha_\ell-\alpha_\ell)^2 {X_{i}^{(\ell)}}^2 \, .
\end{equation}
Moreover
\begin{equation}\label{B2}
|\widehat {\mathbf A}_{j,k}- \widetilde {\mathbf A}_{j,k}| \leq \frac 1 n \sum_{i=1}^n |(X_{i}^{(j)}-\overline{\mathbf{X}}^{(j)}) (X_{i}^{(k)}-\overline{\mathbf{X}}^{(k)}) R_i | \, .
\end{equation}
Combining \eqref{B1} and \eqref{B2}, and using that $|\widehat \alpha_j-\alpha_j|$ converges almost surely to 0, we infer that $\widehat {\mathbf A}_{j,k}- \widetilde {\mathbf A}_{j,k}$ converges almost surely to 0 as soon as 
$$
\limsup_{n \rightarrow \infty} \frac 1 n \sum_{i=1}^n |X_{i}^{(j)}X_{i}^{(k)} X_{i}^{(\ell)} \varepsilon_i | < \infty \quad \text{and}  \quad \limsup_{n \rightarrow \infty} \frac 1 n \sum_{i=1}^n |X_{i}^{(j)}X_{i}^{(k)} {X_{i}^{(\ell)}}^2 | < \infty \quad \text{almost surely},
$$
which is true (applying once again the strong law of large numbers) since  ${\mathbb E}(Y^4)< \infty$ and   ${\mathbb E}({X^{(j)}}^4)< \infty$  for any $j \in \{1, \ldots , p\}$.

Let us now prove that, for  $(j,k) \in\{p+1, \ldots, 2p\}^2$, $\widehat {\mathbf A}_{j,k}- \widetilde {\mathbf A}_{j,k}$ converges almost surely to 0  as $n\rightarrow \infty$. We write
$
\widehat e^{(j-p)}_{i} \widehat e^{(k-p)}_{i} = e^{(j-p)}_{i} e^{(k-p)}_{i} +T_i
$, with
\begin{multline}\label{B3}
|T_i| \leq  |e^{(j-p)}_{i}(\widehat e^{(k-p)}_{i}- e^{(k-p)}_{i})| +|e^{(k-p)}_{i}(\widehat e^{(j-p)}_{i}- e^{(j-p)}_{i})|\\+\frac 1 2 (\widehat e^{(k-p)}_{i}- e^{(k-p)}_{i})^2 +\frac 1 2 (\widehat e^{(j-p)}_{i}- e^{(j-p)}_{i})^2 \, .
\end{multline}
For $ \ell \in\{1, \ldots, p\}$, let $\beta_\ell= {\mathbb E}(X_\ell)-\theta_\ell {\mathbb E}(Y)$ and $ \widehat \beta_\ell= \overline{\mathbf{X}}^{(\ell)}-\widehat \theta_\ell \overline{\mathbf{Y}}$. With these notations
\begin{equation}\label{models}
  X_{k}^{(\ell)}=\beta_\ell + \theta_\ell Y_k + e^{(\ell)}_{k} \quad \text{and}  \quad X_{k}^{(\ell)}= \widehat \beta_\ell + \widehat \theta_\ell Y_k + \widehat e^{(\ell)}_{k}.
\end{equation}
From \eqref{B3} and \eqref{models}, we infer that
\begin{multline}\label{B4}
|T_i| \leq  |e^{(j-p)}_{i}(\widehat \beta_{k-p}- \beta_{k-p})| + |e^{(j-p)}_{i}Y_{i}(\widehat \theta_{k-p}- \theta_{k-p})| \\
+|e^{(k-p)}_{i}(\widehat \beta_{j-p}- \beta_{j-p})| + |e^{(k-p)}_{i}Y_i(\widehat \theta_{j-p}- \theta_{j-p})|\\+ (\widehat \beta_{k-p}- \beta_{k-p})^2 + (\widehat \beta_{j-p}- \beta_{j-p})^2 \\ 
+ (\widehat \theta_{k-p}- \theta_{k-p})^2Y_i^2 + (\widehat \theta_{j-p}- \theta_{j-p})^2 Y_i^2\, .
\end{multline}
Moreover
\begin{equation}\label{B5}
|\widehat {\mathbf A}_{j,k}- \widetilde {\mathbf A}_{j,k}| \leq \frac 1 n \sum_{i=1}^n |(Y_{i}-\overline{\mathbf{Y}})^2 T_i | \, .
\end{equation}
Combining \eqref{B4} and \eqref{B5}, and using that $|\widehat \beta_\ell-\beta_\ell|$ and  $|\widehat \theta_\ell-\theta_\ell|$ converge almost surely to 0, we infer that $\widehat A_{j,k}- \tilde A_{j,k}$ converges almost surely to 0 as soon as, almost surely
$$
\limsup_{n \rightarrow \infty} \frac 1 n \sum_{i=1}^n |Y_i^2 e^{(j-p)}_{i} | < \infty, \
\limsup_{n \rightarrow \infty} \frac 1 n \sum_{i=1}^n |Y_i^3 e^{(j-p)}_{i} | < \infty,  \ \text{and}  \ \limsup_{n \rightarrow \infty} \frac 1 n \sum_{i=1}^n Y_i^4  < \infty,
$$
which is true (applying once again the strong law of large numbers) since  ${\mathbb E}(Y^4)< \infty$ and   ${\mathbb E} ({X^{(j)}}^4 )< \infty$  for any $j \in \{1, \ldots , p\}$.

The fact that, for  $(j,k) \in \{1, \ldots , p\} \times \{p+1, \ldots, 2p\}$, $\widehat {\mathbf A}_{j,k}- \widetilde {\mathbf A}_{j,k}$ converges almost surely to 0  as $n\rightarrow \infty$ may be proved exactly in the same way.
 \qed
 
 \medskip
 
 As an immediate Corollary, we get
 \begin{corollary}\label{CI}
Assume that Assumptions \ref{A_independent} and \ref{A_moment_2} hold, and that ${\mathbb E}(Y^4)< \infty$.
Assume also that $V>0$, where $V$ is defined by \eqref{defV}.  For $\delta \in (0,1)$, let $c_{1-(\delta/2)}$ be the quantile of order $1-\delta$ of the ${\mathcal N}(0,1)$-distribution, and let $\widehat V_n$ be defined by \eqref{hatV}. Then 
$$
\left [  \widehat R^2 -\frac{c_{1-(\delta/2)}\sqrt{\widehat V_n}}{\sqrt n} , \widehat R^2 + \frac{c_{1-(\delta/2)}\sqrt{\widehat V_n}}{\sqrt n} \right] 
$$
is an asymptotic confidence interval for $R^2$ of level $1-\delta$.
\end{corollary}

\begin{Remark}
    As is sensitivity analysis, one can also define 
    $$
    R^2_{T_i}=1-\frac{\mathrm{Var}(\mathrm{Proj}_{V_{(-i)}(Y)}) }{\mathrm{Var}(Y)} \, ,
    $$
    where $V_{(-i)}$ is the linear space of ${\mathbb L}^2$ generated by the variables ${\bf 1},X^{(1)}, \ldots, X^{(i-1)},X^{(i+1)}, \ldots,  X^{(p)}$. Then $R^2_{T_i}$ is the proportion of the variance of $Y$ that is not explained by the best linear predictor based on the variables ${\bf 1},X^{(1)}, \ldots, X^{(i-1)},X^{(i+1)}, \ldots,  X^{(p)}$. Denoting by $R^2_{(-i)}$ the proportion of the variance of $Y$ that is explained by the best linear predictor based on ${\bf 1},X^{(1)}, \ldots, X^{(i-1)},X^{(i+1)}, \ldots,  X^{(p)}$, we immediately see that $R^2_{T_i}=1-R^2_{(-i)}$. Consequently, a confidence interval for $R^2_{T_i}$ is directly obtained from a confidence interval for $R^2_{(-i)}$.
\end{Remark}

To be complete, we discuss the positivity of the variance term $V$.

 \begin{Lemma}\label{V}
 Assume that  $\mathbf M_0$ is invertible, that ${\mathbb E}((X^{(j)}-{\mathbb E}(X^{(j)})^2 \varepsilon^2)< \infty$ and   Assumption \ref{A_moment_3} holds.
 \begin{enumerate}
 \item If $R^2=1$ then $\widehat R^2=1$ almost surely, and $V=0$.
 \item If $R^2=0$ then $V=0$, and $n \widehat R^2$ converges in distribution to $(G_1, \ldots, G_p)(G_{p+1}, \ldots, G_{2p})^t$, where $(G_1, \ldots, G_{2p})$ is the Gaussian random vector with covariance matrix ${\mathbf B}{\mathbf A}{\mathbf B}$.
 \item If $R^2 > 0$ and if  the family 
 \begin{equation}\label{fam}
 ( (Y-{\mathbb E}(Y))^2, ((Y-{\mathbb E}(Y))(X^{(j)}-{\mathbb E}(X^{(j)})))_{1\leq j \leq p} , ((X^{(j)}-{\mathbb E}(X^{(j)}))(X^{(j)}-{\mathbb E}(X^{(j)})))_{1\leq i \leq j \leq p})
 \end{equation}
  is linearly independent, then $V>0$.
 \end{enumerate}
 \end{Lemma}
 
 \begin{Remark}
 The fact that the family 
 \eqref{fam}
  is linearly independent implies that  $(Y, X^{(1)}, \ldots, X^{(p)})$ is also linearly independent,  so that  $R^2 <1$.
 \end{Remark}
 
 \begin{Remark}
 Note that one can easily check that the variable $(G_1, \ldots, G_p)(G_{p+1}, \ldots, G_{2p})^t$ of Item 2 is non negative. Indeed, in that particular case, $ (G_1, \ldots, G_p)^t= \mathbf M_0^{-1} (Z_1, \ldots, Z_p)^t $ and $(G_{p+1}, \ldots, G_{2p})^t= (\mathrm{Var}(Y))^{-1}(Z_1, \ldots, Z_p)^t$, where $(Z_1, \ldots, Z_p)^t$ is a Gaussian vector with covariance matrix $({\mathbf A}_{i,j})_{1\leq i,j\leq p}$.
 The asymptotic distribution of Item 2 may be used to test $H_0: R^2=0$ against $H_1: R^2 \neq 0$, but it is simpler to test the equivalent hypothesis $H_0: \alpha_1= \cdots = \alpha_p=0$ against 
 $H_1: \alpha_i \neq 0$ for some $i \in \{1, \ldots , p\}$ (see Remark \ref{test}). 
 \end{Remark}
 
 \noindent{\bf Proof of Lemma \ref{V}.}
 Item 1 is clear: $R^2=1$ if and only if $Y=\alpha_0 + \alpha_1 X^{(1)} + \cdots + \alpha_p X^{(p)}$ almost surely. In that case $\widehat R^2=1$ almost surely. Consequently $\widehat R^2-R^2=0$ almost surely, and $V=0$.
 We now prove Item 2. $R^2=0$ if and only if $Y-{\mathbb E}(Y)$ is orthogonal to the space generated by $(X^{(1)}-{\mathbb E}(X^{(1)}), \ldots, (X^{(p)}-{\mathbb E}(X^{(p)}))$. This is equivalent to $\alpha_1= \cdots = \alpha_p =0$, which is also equivalent to $\theta_1= \cdots = \theta_p =0$. This implies that $V=0$. Applying Proposition \ref{limitlaw}, we see that
$$
  \sqrt n (\widehat \alpha_1, \ldots, \widehat \alpha_p, \widehat \theta_1, \ldots, \widehat \theta_p)^t  \cvl
 {\mathcal N}_{2p}(0, \B\A\B).
$$
The last assertion of Item 2 follows by noting that $$n \widehat R^2=  \left (\sqrt n (\widehat \theta_2, \ldots, \widehat \theta_p) \right ) \left ( \sqrt n (\widehat \alpha_2, \ldots, \widehat \alpha_p) \right )^t.$$

We now prove Item 3. Since $R^2 > 0$, it follows that the two vectors $(\alpha_1, \ldots, \alpha_p)^t$ and $(\theta_1, \ldots, \theta_p)^t$ are not  equal to $(0, \ldots, 0)^t$. 
Let 
$$
    (a_1, \ldots, a_p, b_1, \ldots, b_p)^t=\B (\theta_1, \ldots, \theta_p, \alpha_1, \ldots, \alpha_p)^t \, .
$$
By definition of the matrix $\B$, and  since $\M_0$ is invertible, we infer that the two vectors $(a_1, \ldots, a_p)^t$ and $(b_1, \ldots, b_p)^t$ are not  equal to $(0, \ldots, 0)^t$. 

Now, by definition, $V=(a_1, \ldots, a_p, b_1, \ldots, b_p)\A(a_1, \ldots, a_p, b_1, \ldots, b_p)^t$. Since $\A$ is the covariance matrix of the vector
$$
\left ((X^{(1)}-{\mathbb E}(X^{(1)}))\varepsilon, \ldots, (X^{(p)}-{\mathbb E}(X^{(p)}))\varepsilon, (Y-{\mathbb E}(Y)) e^{(1)}, \ldots, (Y-{\mathbb E}(Y)) e^{(p)}\right )^t \, ,
$$
we can write
$$
V=\text{Var} \left (a_1(X^{(1)}-{\mathbb E}(X^{(1)}))\varepsilon + \cdots + a_p(X^{(p)}-{\mathbb E}(X^{(p)}))\varepsilon + b_1(Y-{\mathbb E}(Y)) e^{(1)}+ \cdots + b_p(Y-{\mathbb E}(Y)) e^{(p)} \right ) \, .
$$
Hence $V=0$ if and only if 
\begin{equation}\label{lindep}
a_1(X^{(1)}-{\mathbb E}(X^{(1)}))\varepsilon + \cdots + a_p(X^{(p)}-{\mathbb E}(X^{(p)}))\varepsilon = -(b_1(Y-{\mathbb E}(Y)) e^{(1)} + \cdots + b_p(Y-{\mathbb E}(Y)) e^{(p)} ) \quad \text{almost surely.}
\end{equation}
Recall now that
\begin{align*}\label{lindep}
   \varepsilon&= (Y-{\mathbb E}(Y))- \alpha_1(X^{(1)} -{\mathbb E}(X^{(1)})) - \cdots - \alpha_p(X^{(p)} -{\mathbb E}(X^{(p)}))  \\ e^{(k)}&= (X^{(k)} -{\mathbb E}(X^{(k)}))-\theta_k (Y-{\mathbb E}(Y)) \quad  \text{for $k \in\{1, \ldots, p\}$}.
\end{align*}
Without loss of generality, assume that $\alpha_1 \neq 0$. Then, on the left side of \eqref{lindep} the terms 
\begin{multline*}
 -a_1\alpha_1(X^{(1)}-{\mathbb E}(X^{(1)}))^2, -a_2\alpha_1(X^{(1)}-{\mathbb E}(X^{(1)}))(X^{(2)}-{\mathbb E}(X^{(2)})), \\ \ldots, -a_p\alpha_1(X^{(1)}-{\mathbb E}(X^{(1)}))(X^{(p)}-{\mathbb E}(X^{(p)}))
\end{multline*}
appear, but they do not appear in the right side of \eqref{lindep}. Since we assumed that the family \eqref{fam} is linearly independent, we infer that $V=0$ implies $a_1= \cdots = a_p =0$. Since we know that $(a_1, \ldots, a_p)^t$ is not  equal to $(0, \ldots, 0)^t$, we conclude that $V>0$.  \qed

\section{Simulations} \label{Sec4}
\setcounter{equation}{0}
In this section, we evaluate the performance of the confidence interval of Corollary \ref{CI} on different sets of simulations, and we compare this performance with that of the confidence interval that assumes that the vector $(Y, X^{(1)}, \ldots X^{(p)})$ is Gaussian.

\medskip

We shall consider different models. For each model, we shall estimate the coverage level of two confidence intervals of level $95\%$ for the $R^2$, for $n=200$ to $n=1000$, via  a basic Monte-Carlo procedure (with $N=3000$ repetitions). 

The first confidence interval (CI1) is the one described in Corollary \ref{CI} (with the small change that we take the quantile $t_{n, 0.975}$ of the Student distribution St($n$) instead of the quantile $c_{0.975}$; this has no theoretical justification, but can improve a bit the coverage level for small $n$). The second confidence interval (CI2) is the non-asymptotic confidence interval ci.R2 of the R package MBESS, as described by  Kelley \cite{K}, and based on a precise approximation of the distribution of  
$\widehat R^2$ given by Lee \cite{L1}, \cite{L2}, when the vector $(Y,X^{(1)}, \dots, X^{(p)})$ is Gaussian. It is clearly indicated in Kelley \cite{K2} (Discussion, page 552-553) that CI2 is not robust to the non-normality of $(Y,X^{(1)}, \dots, X^{(p)})$, which will be confirmed by the simulations.

\subsection{An example where $(Y,X^{(1)},X^{(2)})$ is Gaussian}
We consider here the model
\begin{equation}\label{modG}
Y=0.5+ 0.5 X^{(1)} + X^{(2)} + \varepsilon \, ,
\end{equation}
where $X_1, X_2, \varepsilon$ are i.i.d. with ${\mathcal N}(0,1)$ distribution. Let $V$ be the sub-space of ${\mathbb L}^2$ generated by ${\bf 1}, X^{(1)}, X^{(2)}$, so that $\mbox{Proj}_V(Y)=0.5+ 0.5 X^{(1)} + X^{(2)}$. One can then easily check that 
$R^2=5/9$. 

The estimated coverage levels for CI1 and CI2 are given in Table \ref{tab1} below.

 \begin{table}
\center
\begin{tabular}{|c|c|c|c|c|c|c|c|c|c|}
\hline
 $n $  &    200  & 300 & 400 & 500 & 600& 700  &800 & 900 &1000 \\
\hline 
CI1  & 0.942  &  0.945 & 0.952 &  0.949 & 0.951 & 0.947 &0.947 & 0.952& 0.95\\
\hline
CI2 &0.949 &0.952 & 0.954 & 0.95 & 0.953 &0.951&0.949 &0.951 &0.951  \\
\hline
\end{tabular}
\caption{Estimated coverage levels of CI1 and CI2 at level $95\%$ for model  \eqref{modG} with 
$\varepsilon \sim  {\mathcal N}(0,1)$ and $\mbox{Proj}_V(Y)=0.5+  0.5X^{(1)}+  X^{(2)}$.}
\label{tab1}
\end{table}

We see that the estimated coverage level of CI2 is always close to 0.95, which is not a surprise since CI2 is based on a precise approximation of the distribution of  $\widehat R^2$ in the case where $(Y,X^{(1)},X^{(2)})$ is Gaussian. We see that CI1 also gives very good results, with estimated coverage levels between $0.945$ and $0.952$ as soon as $n\geq 300$.

\subsection{An example where the error term has a Student distribution}
We consider here a slight modification of  the model \eqref{modG}, where $\varepsilon \sim$ St(10). 
Again,  $V$ is the sub-space of ${\mathbb L}^2$ generated by ${\bf 1}, X^{(1)}, X^{(2)}$, so that $\mbox{Proj}_V(Y)=0.5+ 0.5 X^{(1)} + X^{(2)}$. One can then easily check that 
$R^2=0.5$. 

The estimated coverage levels for CI1 and CI2 are given in Table \ref{tab2} below.

 \begin{table}
\center
\begin{tabular}{|c|c|c|c|c|c|c|c|c|c|}
\hline
 $n $  &    200  & 300 & 400 & 500 & 600& 700  &800 & 900 &1000 \\
\hline 
CI1  & 0.934  &  0.937 & 0.941 &  0.942 & 0.947 & 0.949 &0.95 & 0.945& 0.95\\
\hline
CI2 &0.936 &0.932 & 0.933 & 0.936 & 0.937 &0.936 &0.939 &0.935 &0.934  \\
\hline
\end{tabular}
\caption{Estimated coverage levels of CI1 and CI2 at level $95\%$ for model  \eqref{modG} with 
$\varepsilon \sim$ St(10) and  $\mbox{Proj}_V(Y)=0.5+  0.5X^{(1)}+  X^{(2)}$.}
\label{tab2}
\end{table}

We see that the estimated coverage level of CI2 is always  between 0.932 and 0.94, and is not getting closer to 0.95 as $n$ increases. This confirms that CI2 is not robust to non normality, even when the distribution of the error is symmetric.  We see that CI1 is always better than CI2 as soon as $n\geq 300$, with a coverage level between $0.94$ and $0.95$ when $n\geq 400$.  

\subsection{An heteroscedastic example}
We continue with a modification of the model \eqref{modG}. We consider the model 
\begin{equation}\label{modGH}
Y=0.5+ 0.5 X^{(1)} + X^{(2)} + \left(\sqrt{0.2+0.8{X^{(1)}}^2}\right) e \, ,
\end{equation}
where $X^{(1)}, X^{(2)}, e$ are i.i.d. with ${\mathcal N}(0,1)$ distribution. Let $V$ be the sub-space of ${\mathbb L}^2$ generated by ${\bf 1}, X^{(1)}, X^{(2)}$, so that $\mbox{Proj}_V(Y)=0.5+ 0.5 X^{(1)} + X^{(2)}$. Again, one can easily check that 
$R^2=5/9$. 

The estimated coverage levels for CI1 and CI2 are given in Table \ref{tab9} below.

 \begin{table}
\center
\begin{tabular}{|c|c|c|c|c|c|c|c|c|c|}
\hline
 $n $  &    200  & 300 & 400 & 500 & 600& 700  &800 & 900 &1000 \\
\hline 
CI1  & 0.91  &  0.926 &  0.928 & 0.935 & 0.937 & 0.939 & 0.944 &0.944 & 0.946\\
\hline
CI2 &0.878 &0.878 & 0.878&0.877 & 0.877 & 0.887 &0.884&0.885 &0.887   \\
\hline
\end{tabular}
\caption{Estimated coverage levels of CI1 and CI2 at level $95\%$ for the heteroscedastic model  \eqref{modGH} 
with $\mbox{Proj}_V(Y)=0.5+  0.5X^{(1)}+  X^{(2)}$.}
\label{tab9}
\end{table}

We see that the estimated coverage level of CI2 is always around 0.88. Again, this confirms that CI2 is not robust to non normality.  We see that CI1 is always better than CI2, with a coverage level greater than $0.935$ for $n\geq 500$, and greater than $0.94$ for $n\geq 800$.

\subsection{An example where the model is misspecified}
 We consider the model 
\begin{equation}\label{modMiss}
Y=X^2+ e \, ,
\end{equation}
where $X \sim {\mathcal N}(0,1)$, $e \sim {\mathcal N}(0,1)$  and $e$ is independent of $X$. Let $V$ be the sub-space of ${\mathbb L}^2$ generated by ${\bf 1}, X, |X|$. Elementary computations show that $\mbox{Proj}_V(Y)= \alpha_0 + \alpha_1 X + \alpha_2 |X|$ with $\alpha_0\simeq  -0.752, \alpha_1=0, \alpha_2 \simeq 2.196$,  and $R^2 \simeq 0.584$. The model is misspecified in the sense that $X^2$ does not belong to $V$. It follows that the error term $\varepsilon= Y- \mbox{Proj}_V(Y)$ is such that
$$
{\mathbb E}(\varepsilon | X)=X^2-\alpha_0- \alpha_2 |X| \neq 0 \, .
$$ 

The estimated coverage levels for CI1 and CI2 are given in Table \ref{tab10} below.

 \begin{table}
\center
\begin{tabular}{|c|c|c|c|c|c|c|c|c|c|}
\hline
 $n $  &    200  & 300 & 400 & 500 & 600& 700  &800 & 900 &1000 \\
\hline 
CI1  & 0.935  &  0.941 &  0.943 & 0.943 & 0.943 & 0.948 & 0.947 &0.945 & 0.947\\
\hline
CI2 &0.875 &0.871 & 0.885&0.868 & 0.861 & 0.87 &0.866&0.874 &0.878   \\
\hline
\end{tabular}
\caption{Estimated coverage levels of CI1 and CI2 at level $95\%$ for  model  \eqref{modMiss} 
with $\mbox{Proj}_V(Y)\simeq  -0.752 + 2.196|X|$.}
\label{tab10}
\end{table}

We see that the estimated coverage level of CI2 is always around 0.87. Again, this confirms that CI2 is not robust to non normality.  We see that CI1 is always better than CI2, with a coverage level greater than $0.94$ for $n\geq 300$.

\subsection{An example of polynomial regression}
 We consider again the model 
\eqref{modMiss}, but now $V$ is the sub-space of ${\mathbb L}^2$ generated by ${\bf 1}, X, X^2$. It is then obvious that  $\mbox{Proj}_V(Y)= \alpha_0 + \alpha_1 X + \alpha_2 X^2$ with $\alpha_0=0, \alpha_1=0, \alpha_2 =1$,  and $R^2 =2/3$. 

The estimated coverage levels for CI1 and CI2 are given in Table \ref{tab11} below.

 \begin{table}
\center
\begin{tabular}{|c|c|c|c|c|c|c|c|c|c|}
\hline
 $n $  &    200  & 300 & 400 & 500 & 600& 700  &800 & 900 &1000 \\
\hline 
CI1  & 0.914  &  0.926 &  0.927 & 0.933 & 0.934 & 0.939 & 0.94 &0.94 & 0.942\\
\hline
CI2 &0.75 &0.755 & 0.752 &0.751 & 0.754 & 0.759 &0.74&0.755 &0.746   \\
\hline
\end{tabular}
\caption{Estimated coverage levels of CI1 and CI2 at level $95\%$ for  model  \eqref{modMiss} 
with $\mbox{Proj}_V(Y)=X^2$.}
\label{tab11}
\end{table}

We see that the estimated coverage level of CI2 is always around 0.75. Again, this confirms that CI2 is not robust to non normality.  We see that CI1 is always much better than CI2, with a coverage level greater than $0.93$ for $n\geq 500$ and greater than $0.94$ for $n\geq 800$.

\subsection{An example of Poisson regression}
In this example, $V$ is the subspace of ${\mathbb L}^2$ generated by ${\bf 1}, X^{(1)}, X^{(2)}$, where $X^{(1)}, X^{(2)}$ are independent, $X^{(1)}$ is uniformly distributed over $[0,1]$, and $X^{(2)}$ is exponentially distributed with parameter 1. The response variable $Y$ is a count variable, whose conditional distribution given $X^{(1)}, X^{(2)}$ is a Poisson distribution with parameter 
$0.5+X^{(1)}+X^{(2)}$. In that case, we have 
$\mbox{Proj}_V(Y)= \mathbb{E}(Y|X^{(1)}, X^{(2)})= \alpha_0 + \alpha_1 X^{(1)} + \alpha_2 X^{(2)}$ with $\alpha_0=0.5, \alpha_1=2, \alpha_2 =1$,  and $R^2 =13/37$. Note that this is again an heteroscedastic model, since 
$$
\mathrm{Var}(\varepsilon|X^{(1)}, X^{(2)})= \mathrm{Var}(Y-\mbox{Proj}_V(Y)|X^{(1)}, X^{(2)})= 0.5+X^{(1)}+X^{(2)}\, .
$$

The estimated coverage levels for CI1 and CI2 are given in Table \ref{tab12} below.

\begin{table}
\center
\begin{tabular}{|c|c|c|c|c|c|c|c|c|c|}
\hline
 $n $  &    200  & 300 & 400 & 500 & 600& 700  &800 & 900 &1000 \\
\hline 
CI1  & 0.922  &  0.932 &  0.933 & 0.936 & 0.94 & 0.941 & 0.946 &0.943 & 0.943\\
\hline
CI2 &0.872 &0.863 & 0.857 &0.861 & 0.863 & 0.853 &0.862 &0.858 &0.861   \\
\hline
\end{tabular}
\caption{Estimated coverage levels of CI1 and CI2 at level $95\%$ for The Poisson regression  model  
with $\mbox{Proj}_V(Y)=0.5+ X^{(1)} +X^{(2)}$.}
\label{tab12}
\end{table}

We see that the estimated coverage level of CI2 is always around 0.86. Again, this confirms that CI2 is not robust to non normality.  We see that CI1 is always much better than CI2, with a coverage level greater than $0.93$ for $n\geq 300$ and greater than $0.94$ for $n\geq 600$.

\section{Partial correlations, vector-valued random variables,  dependent sequences} \label{Sec5}

In this section, we show that the approach developed in Section 3 can be adapted, with minor changes, to the case of partial correlation coefficients, to the case of vector-valued random variables, and to the case where the observations come from a stationary $\alpha$-mixing sequence.

\subsection{Partial correlations} \label{Sec6.1}
In this section, we show that the approach of Section 3 works also for  partial correlations. For the sake of simplicity, we shall only consider the case where there are two real-valued random variables $X,Y$ in ${\mathbb L}^2$, and a vector of confounding random variables ${\bf Z}=(Z^{(1)}, \ldots, Z^{(p)})$ (each $Z^{(k)}$ being square integrable). We assume that the variables ${\bf 1}, X,Y,Z^{(1)}, \ldots, Z^{(p)}$ are linearly independent (as a family of random variables in ${\mathbb L}^2$).

Let $V({\bf Z})$ be the subspace of ${\mathbb L}^2$ generated by ${\bf 1}, Z^{(1)}, \ldots , Z^{(p)}$, and let
\begin{align*}
X({\bf Z})& =X- \mbox{Proj}_{V({\bf Z})}(X) = X-a_0 -a_1 Z^{(1)} - \cdots - a_p Z^{(p)}\, ,\\
Y({\bf Z})& =Y- \mbox{Proj}_{V({\bf Z})}(Y) = Y-b_0 -b_1 Z^{(1)} - \cdots - b_p Z^{(p)}\, .
\end{align*}

The square of the partial correlation between $X$ and $Y$ with respect to ${\bf Z}$ is then simply the $R^2$ coefficient between $X({\bf Z})$ and $Y({\bf Z})$, which we denote by $R^2(X({\bf Z}), Y({\bf Z}))$. It is a way tho measure the square of the correlation between $X$ and $Y$ after removing the correlation due the the confounding vector ${\bf Z}$.

Let then 
\begin{align*}
X({\bf Z})&=\alpha_0 + \alpha Y({\bf Z})  + \varepsilon^{(1)} \, , \\
Y({\bf Z})&=\theta_0 + \theta X({\bf Z}) + \varepsilon^{(2)} \, ,
\end{align*}
 where $\varepsilon^{(1)}$ is orthogonal to ${\bf 1}, Y({\bf Z})$, and $\varepsilon^{(2)}$ is orthogonal to ${\bf 1}, X({\bf Z})$. From Lemma \ref{scalarproduct}, we know that 
 $$
 R^2(X({\bf Z}), Y({\bf Z}))=\alpha \theta \, .
 $$

 Let now $(X_i, Y_i, {\bf Z}_i)_{1\leq i \leq n}$ be $n$ independent copies of $(X,Y, {\bf Z})$.  
 The empirical estimator of the coefficient $R^2(X({\bf Z}), Y({\bf Z}))$ is then
 $$
 \widehat R^2(X({\bf Z}), Y({\bf Z})) = \widehat \alpha \widehat \theta \, ,
 $$
 where 
 \begin{align*}
 \widehat \alpha &= \frac{\sum_{i=1}^n (X_i- \hat a_0 - \cdots -\hat a_p Z_i^{(p)})(X_i- \hat b_0 - \cdots -\hat b_p Z_i^{(p)})}
 {\sum_{i=1}^n (X_i- \hat a_0 - \cdots -\hat a_p Z_i^{(p)})^2} \\
 \widehat \theta &= \frac{\sum_{i=1}^n (X_i- \hat a_0 - \cdots -\hat a_p Z_i^{(p)})(X_i- \hat b_0 - \cdots -\hat b_p Z_i^{(p)})}
 {\sum_{i=1}^n (Y_i- \hat b_0 - \cdots -\hat b_p Z_i^{(p)})^2} \, .
 \end{align*}
 To describe the asymptotic behavior of $\sqrt n \left(\widehat R^2(X({\bf Z}), Y({\bf Z}))-  R^2(X({\bf Z}), Y({\bf Z})\right )$, we need to define the two $2\times 2$ symmetric matrices $\bf A$ and $\bf B$ as follows : 

 \begin{itemize}
 
 \item ${\bf A}_{1,1}= {\mathbb E}\left( X({\bf Z})^2{\varepsilon^{(1)}}^2\right )$, ${\bf A}_{2,2}= {\mathbb E}\left( Y({\bf Z})^2{\varepsilon^{(2)}}^2\right )$, ${\bf A}_{1,2}= {\mathbb E}\left( X({\bf Z})Y({\bf Z}){\varepsilon^{(1)}} {\varepsilon^{(2)}}\right )$.

 \item ${\bf B}_{1,1}= (\text{Var}\left (X({\bf Z})) \right)^{-1} $, ${\bf B}_{2,2}=(\text{Var}(Y({\bf Z})))^{-1} $, ${\bf B}_{1,2}= 0$.
 \end{itemize}

 \begin{Proposition}\label{limitlawR2Par}
Assume  that the variables ${\bf 1}, X,Y,Z^{(1)}, \ldots, Z^{(p)}$ are linearly independent. Assume  that ${\mathbb E}\left({Z^{(j)}}^2 X({\bf Z})^2\right )<\infty$, ${\mathbb E}\left({Z^{(j)}}^2 Y({\bf Z})^2\right )<\infty$ for any $j \in \{1, \ldots , p\}$, and  that ${\mathbb E}\left( X({\bf Z})^2{\varepsilon^{(1)}}^2\right )<\infty$,  ${\mathbb E}\left( Y({\bf Z})^2{\varepsilon^{(2)}}^2\right )< \infty$.
Then
$$
   \sqrt n \left(\widehat R^2(X({\bf Z}), Y({\bf Z}))-  R^2(X({\bf Z}), Y({\bf Z})\right )  \cvl {\mathcal N}(0, V).
$$
where 
\begin{equation*}\label{defV2}
V= (  \theta, \alpha){\mathbf B}{\mathbf A}{\mathbf B}(  \theta, \alpha )^t \, .
\end{equation*}
\end{Proposition}

\begin{Remark}
The matrices $\A$ and $\B$ can be consistently estimated by their empirical counterparts based on the residuals, provided 
${\mathbb E}(X^4)< \infty$, ${\mathbb E}(Y^4)< \infty$ and ${\mathbb E}({Z^{(j)}}^4)< \infty$ for any $j \in \{1, \ldots , p\}$ (as done in the proof Proposition \ref{estim}).
As in Proposition \ref{estim} this provides a consistent estimator
of the limiting variance $V$.
\end{Remark}

\noindent{\bf Proof of Proposition \ref{limitlawR2Par}.}
Let $\widetilde R^2(X({\bf Z}), Y({\bf Z}))=\widetilde \alpha \widetilde \theta$, where 
\begin{align*}
 \widetilde \alpha &= \frac{\sum_{i=1}^n (X_i- a_0 - \cdots -a_p Z_i^{(p)})(X_i-  b_0 - \cdots - b_p Z_i^{(p)})}
 {\sum_{i=1}^n (X_i-  a_0 - \cdots - a_p Z_i^{(p)})^2} \\
 \widetilde \theta &= \frac{\sum_{i=1}^n (X_i-  a_0 - \cdots - a_p Z_i^{(p)})(X_i-  b_0 - \cdots - b_p Z_i^{(p)})}
 {\sum_{i=1}^n (Y_i-  b_0 - \cdots - b_p Z_i^{(p)})^2} \, .
 \end{align*}
 Applying Proposition \ref{limitlawR2}, we see that Proposition \ref{limitlawR2Par} is true with $\widetilde R^2(X({\bf Z}), Y({\bf Z}))$ instead of the empirical estimator $\widehat R^2(X({\bf Z}), Y({\bf Z}))$. To conclude, it suffices to prove that 
 $$
 \sqrt n \left ((\widehat \alpha, \widehat \theta)- (\widetilde \alpha, \widetilde \theta) \right ) \quad \text{converges in probability to 0.}
 $$
 Clearly, it suffice to prove that the two quantities 
 \begin{multline*}
 \sqrt{n} \Big ( \frac 1 n \sum_{i=1}^n (X_i- a_0 - \cdots -a_p Z_i^{(p)})(X_i-  b_0 - \cdots - b_p Z_i^{(p)})\\- \frac 1 n \sum_{i=1}^n (X_i- \hat a_0 - \cdots -\hat a_p Z_i^{(p)})(X_i-  \hat b_0 - \cdots - \hat b_p Z_i^{(p)}) \Big )
 \end{multline*}
 and 
 $$
 \sqrt{n} \left ( \frac 1 n \sum_{i=1}^n (X_i- a_0 - \cdots -a_p Z_i^{(p)})^2\\- \frac 1 n \sum_{i=1}^n (X_i- \hat a_0 - \cdots -\hat a_p Z_i^{(p)})^2 \right )
 $$
converge in probability to 0. We prove the second convergence, the proof of the first one being similar. By Pythagoras theorem in ${\mathbb R}^n$, 
\begin{multline*}
\left |\frac 1 {\sqrt n} \sum_{i=1}^n (X_i- a_0 - \cdots -a_p Z_i^{(p)})^2- \frac 1 {\sqrt n} \sum_{i=1}^n (X_i- \hat a_0 - \cdots -\hat a_p Z_i^{(p)})^2 \right | \\ \leq p \left(
\sum_{j=1}^p \sqrt n (\hat a_j -a_j)^2 \left (\frac 1 n\sum_{i=1}^n {Z_i^{(j)}}^2 \right )\right) \, .
\end{multline*}
The result follows, since by Proposition \ref{cltP} $\sqrt n (\hat a_i -a_i)^2$ converges in probability to $0$. \qed

\subsection{Vector-valued random variables}
\label{sec:ext}

\setcounter{equation}{0}

In this section, we assume that $(Y, X^{(1)}, \ldots X^{(p)})$ is a random vector with value in ${\mathbb H}^{p+1}$, where ${\mathbb H}$ is a real separable Hilbert space with inner product $\langle \cdot , \cdot \rangle$ and norm $\|\cdot \|_{\mathbb H}$.
Let ${\mathbb L}^2({\mathbb H})$ be the space of ${\mathbb H}$-valued random variables $Z$ such that ${\mathbb E}(\|Z \|^2_{\mathbb H}) < \infty$, and recall that ${\mathbb L}^2({\mathbb H})$ is also a separable Hilbert space with inner product ${\mathbb E}(\langle \cdot , \cdot \rangle)$.

We assume in the following that all the variables $Y, X^{(1)}, \ldots X^{(p)}$ belong to ${\mathbb L}^2({\mathbb H})$.
Let $V$ be the subspace of  ${\mathbb L}^2({\mathbb H})$ generated by the variables $(X^{(1)}, \ldots, X^{(p)})$, and let  $\mbox{Proj}_V(Y)$ be the orthogonal projection of $Y$ onto $V$. 
To measure the quality of approximation of $Y$ by $\text{Proj}_{V}(Y)$, one can consider the coefficient $R^2$ defined by (assuming that ${\mathbb E}(\|Y \|^2_{\mathbb H})>0$):
$$
R^2= \frac{{\mathbb E}\left ( \|\text{Proj}_{V}(Y)\|^2_{\mathbb H}\right )}{{\mathbb E} \left ( \|Y\|^2_{\mathbb H}\right )} \, .
$$

Let
$(Y_i, X_{i}^{(1)}, \ldots X_{i}^{(p)})_{1 \leq i \leq n}$ be $n$ independent copies of the vector $(Y, X^{(1)}, \ldots X^{(p)})$, and 
let $V_i$ be the ${\mathbb L}^2({\mathbb H})$-subspace generated by the variables  $(X^{(1)}_i, \ldots, X^{(p)}_i)$. Let $\varepsilon=Y- \mbox{Proj}_{V}(Y)$ and $ \varepsilon_i= Y_i- \mbox{Proj}_{V_i}(Y_i)$.
 We can then write
$$
Y_i= \mbox{Proj}_{V_i}(Y_i)+ \varepsilon_i.
$$
where the variables $\varepsilon_i$ are such that ${\mathbb E}(\langle X_{i}^{(j)}, \varepsilon_i\rangle )=0$ for any $j \in \{1, \ldots , p\}$.

Assume that  $ X^{(1)}, \ldots,  X^{(p)}$ are linearly independent (as a family of random variables in  ${\mathbb L}^2({\mathbb H})$).
Then one can uniquely write
$$
Y_i= \alpha_1 X_{i}^{(1)}+ \cdots + \alpha_p X_{i}^{(p)} +\varepsilon_i  \, .
$$
Let $\bm{\alpha}=(\alpha_1,\ldots,\alpha_p)^t\in \mathbb{R}^{p}$, and $\mathbf M$ be the matrix such that  \begin{equation}\label{defMbis}
\mathbf M_{j,\ell}={\mathbb E}(\langle X^{(j)}, X^{(\ell)}\rangle )\, .
\end{equation}
Since $ X^{(1)}, \ldots,  X^{(p)}$ are linearly independent, $\mathbf M$ is invertible and one can easily show that
$$
\bm{\alpha}= \mathbf M^{-1}({\mathbb E}(\langle X^{(1)}, Y\rangle) \ldots, {\mathbb E}(\langle X^{(p)}, Y\rangle ))^t \, .
$$
Let now $\widehat{\mathbf M}$ be the matrix such that  $\widehat{\mathbf M}_{j,\ell}=n^{-1}\sum_{i=1}^n\langle X_i^{(j)}, X_i^{(\ell)}\rangle $. By the strong law of large numbers
\begin{equation}\label{lgn1bis}
\lim _{n \rightarrow \infty} \widehat{\mathbf M} = \mathbf M \ \text{almost surely.}
\end{equation}
Since  $\mathbf M$  is invertible,  for $n$ large enough, $\widehat{\mathbf M}$ is invertible. Consequently,  one can define  
$$
\widehat{\bm{\alpha}}=(\widehat \alpha_1, \ldots, \widehat \alpha_p)^t= \widehat{\mathbf M}^{-1}\left (\frac 1 n \sum_{i=1}^n \langle X^{(1)}_i, Y_i\rangle, \ldots, \frac 1 n \sum_{i=1}^n \langle X^{(p)}_i, Y_i\rangle \right )^t \, .
$$
It follows that 
\begin{equation}\label{estbis}
\widehat{\bm{\alpha}}-\bm{\alpha}=(\widehat \alpha_1, \ldots, \widehat \alpha_p)^t-(\alpha_1, \ldots, \alpha_p)^t = \widehat{\mathbf M}^{-1}\left (\frac 1 n \sum_{i=1}^n \langle X^{(1)}_i, \varepsilon_i\rangle, \ldots, \frac 1 n \sum_{i=1}^n \langle X^{(p)}_i, \varepsilon_i\rangle \right )^t \, .
\end{equation}
Then, proceeding as in Section \ref{sec2.2}, we see that 
$\widehat{\bm{\alpha}}$ converges almost surely to $\bm{\alpha}$, and that: if for any $j \in \{0, \ldots , p\}$, one has  ${\mathbb E}(\langle {X^{(j)}},  \varepsilon \rangle^2)< \infty$,
 then
$$
\sqrt n (\widehat{\bm{\alpha}}-\bm{\alpha})
\cvl {\mathcal N}_{p}(0, \mathbf M^{-1}\mathbf M(\varepsilon)\mathbf M^{-1})\, ,
$$
where $\mathbf M({\varepsilon})$ is the matrix defined by   
$ \mathbf M(\varepsilon)_{j,\ell}={\mathbb E}(\langle X^{(j)}, \varepsilon \rangle \langle X^{(\ell)},  \varepsilon \rangle )$ for any  $0\leq j,\ell\leq p$. 

One can then estimate the matrix $\mathbf M(\varepsilon)$ as in Section \ref{sec:est}.
Define the residuals $\widehat \varepsilon_i$ by:
$$
\widehat \varepsilon_i
=Y_i-\widehat \alpha_1 X_{i}^{(1)} - \cdots - \widehat \alpha_p X_{i}^{(p)}. 
$$ 
Let then
$$
\widehat {\mathbf M(\varepsilon)}_{j,\ell}= \frac 1 n \sum_{i=1}^n \langle X_{i}^{(j)}, \widehat \varepsilon_i \rangle \langle X_{i}^{(\ell)},  \widehat \varepsilon_i \rangle  \, .
$$
As in Proposition \ref{consist}, the following consistency result holds:
assume that   ${\mathbb E}(\|X^{(j)}\|_{\mathbb H}^2 \|\varepsilon \|_{\mathbb H}^2)< \infty$ and ${\mathbb E}(\|X^{(j)}\|_{\mathbb H}^4)< \infty$ for any  $j \in \{0, \ldots , p\}$.
Then, for any $j,\ell \in \{0, \ldots , p\}$,
$$
\lim_{n \rightarrow \infty} \widehat {\mathbf M(\varepsilon)}_{j,\ell}=\mathbf M(\varepsilon)_{j,\ell} \quad \text{almost surely}.
$$
Consequently, since  $\mathbf M$ is invertible, 
$$
\lim_{n \rightarrow \infty} \widehat{\mathbf M}^{-1} \widehat{\mathbf M(\varepsilon)} \widehat{\mathbf M}^{-1}= {\mathbf M}^{-1} {\mathbf M(\varepsilon)} {\mathbf M}^{-1} \quad \text{almost surely}.
$$
The empirical estimator of $R^2$ is 
$$
\widehat R^2= \frac{\sum_{i=1}^n \|\mbox{Proj}_{V_i}(Y_i)\|^2_{\mathbb H}}{\sum_{i=1}^n \|Y_i\|^2_{\mathbb H}} \, .
$$
 Let $ \bm{\theta} =(\theta_1, \ldots, \theta_p)^t$, with
$\widehat{\bm{\theta}} =(\widehat \theta_1, \ldots, \widehat \theta_p)^t$
$$
\theta_j= \frac{{\mathbb E}(\langle X^{(j)}, Y\rangle)}{{\mathbb E}(\|Y\|^2_{\mathbb H})} \quad \text{and} \quad \widehat \theta_j = \frac{\sum_{i=1}^n \langle X^{(j)}_i, Y_i\rangle}{\sum_{i=1}^n\|Y_i\|^2_{\mathbb H}} \quad \text{for $j \in \{1, \ldots , p\}$}.
$$
Then, as in Lemma \ref{scalarproduct} one can prove that
$
 R^2= {\bm{\theta}}^t{\bm{\alpha}}
$
and
$
\widehat R^2= {\widehat{\bm{\theta}}}^t{\widehat{\bm{\alpha}}}
$.
Let 
\begin{equation}\label{def:ekbis}
e^{(k)} = X^{(k)}-\mbox{Proj}_{W}(X^{(k)}) \quad \text{and} \quad e_{i}^{(k)} = X_{i}^{(k)}-\mbox{Proj}_{W_i}(X_{i}^{(k)}),
\end{equation}
where $W$ is the subspace of ${\mathbb L}^2({\mathbb H})$ generated by $Y$, and $W_i$ is the subspace of ${\mathbb L}^2({\mathbb H})$ generated by $Y_i$. Let also ${\bm e}^{(k)}=(e_{1}^{(k)}, \ldots, e_{n}^{(k)})^t $.
As in Proposition 
\ref{limitlawR2}, one can prove the following result:
if  ${\mathbb E}(\langle X^{(j)}, \varepsilon\rangle^2)< \infty$  and  ${\mathbb E}(\langle Y, e^{(j)}\rangle ^2)< \infty$ for all $j \in \{1  , \ldots , p\}$. Then
$$
  \sqrt n \left ( \widehat R^2 - R^2\right )  \cvl {\mathcal N}(0, V).
$$
where 
\begin{equation*}\label{defVbis}
V= (  \theta_1, \ldots,  \theta_p,  \alpha_1, \ldots,  \alpha_p){\mathbf B}{\mathbf A}{\mathbf B}(  \theta_1, \ldots,  \theta_p,  \alpha_1, \ldots, \alpha_p)^t \, .
\end{equation*}
Here, the matrices ${\mathbf A}$ and ${\mathbf B}$ are as follows:
\begin{itemize}
\item If $(j,k) \in \{1, \ldots, p\}^2$ then 
$$\A_{j,k}={\mathbb E}(\langle X^{(j)}, \varepsilon\rangle \langle X^{(k)}, \varepsilon\rangle);$$ 
\item If $(j,k) \in \{p+1, \ldots, 2p\}^2$ then $$\A_{j,k}={\mathbb E}(\langle Y,  e^{(j-p)}\rangle  \langle Y,  e^{(k-p)}\rangle);$$ \item  $(j,k) \in \{1, \ldots , p\} \times \{p+1, \ldots, 2p\}$ then $$\A_{j,k}={\mathbb E}(\langle X^{(j)},\varepsilon\rangle  \langle Y,  e^{(k-p)}\rangle ).$$
\end{itemize}
Let $\delta_{j,k}=0$ if $j\not=k$ and $\delta_{j,j}=1$ and let $\B$ be the $2p\times 2p$ matrix defined as follows: 
\begin{itemize}
    \item 
If $(j,k) \in \{1, \ldots, p\}^2$ then 
$\B_{j,k}=(\mathbf M^{-1})_{j,k}$ (see  \eqref{defMbis} for the definition of $\mathbf M$);
\item If $(j,k) \in \{p+1, \ldots, 2p\}^2$ then $\B_{j,k}=({\mathbb E}(\|Y\|_{\mathbb H}^2))^{-1}\delta_{j,k}$; 
\item If $(j,k) \in \{1, \ldots , p\} \times \{p+1, \ldots, 2p\}$ then $\B_{j,k}=0$. 
\end{itemize}
Now the matrices $\A$ and $\B$ can be consistently estimated by their empirical counterparts based on the residuals, provided ${\mathbb E}(\|Y\|_{\mathbb H}^4)< \infty$ and ${\mathbb E}(\|X^{(j)}\|_{\mathbb H}^4)< \infty$ for any $j \in \{1, \ldots , p\}$ (as done in the proof Proposition \ref{estim}).
As in Proposition \ref{estim} this provides a consistent estimator
of the limiting variance $V$.

\subsection{Dependent sequences}
In this section, we consider the case where the observations come from a strictly stationary ergodic sequence, and give conditions under which Proposition \ref{limitlaw} remains valid (with a different limiting covariance matrix). This will immediately imply the asymptotic normality of $\sqrt n  ( \widehat R^2 -R^2 )$ as in Proposition \ref{limitlawR2}. For simplicity, we consider only the case of mixing sequences in the sense of Rosenblatt  \cite{Ros56}, but similar results under other conditions can of course be obtained.

Let $(Y_i, X_{i}^{(1)}, \ldots ,X_{i}^{(p)})_{1 \leq i \leq n}$ be $n$ ${\mathbb R}^{p+1}$-valued random vectors obtained from a strictly stationary sequence 
$(Y_i, X_{i}^{(1)}, \ldots ,X_{i}^{(p)})_{i \in {\mathbb Z}}$.
Let ${\mathcal F}_0= \sigma(Y_i, X_{i}^{(1)}, \ldots ,X_{i}^{(p)}, i \leq 0)$ and ${\mathcal F}_k= \sigma(Y_k, X_{k}^{(1)}, \ldots ,X_{k}^{(p)})$. We denote by $\alpha(k)$ the $\alpha$-mixing coefficient between the $\sigma$-algebras ${\mathcal F}_0$ and ${\mathcal F}_k$:
\begin{equation}\label{defalpha}
\alpha(k) = \alpha({\mathcal F}_0, {\mathcal F}_k)=\sup_{A \in {\mathcal F}_0, B \in {\mathcal F}_k} \left | {\mathbb P}(A\cap B)- {\mathbb P}(A) {\mathbb P}(B)\right | \, .
\end{equation}
For any real-valued random variables $Z$, let $Q_Z$ be the inverse cadlag of the tail function $t \rightarrow {\mathbb P}(|Z|>t)$ (note that, by definition, $Q_Z$ is non increasing on $[0,1]$).

Keeping the notations of Section \ref{Sec3}, we replace Assumption \ref{A_moment_3} by: for all $j \in \{1, \ldots , p\}$

\begin{equation}\label{condalpha}
    \sum_{k\geq 0} \int_0^{\alpha(k)}Q^2_{(X^{(j)}-{\mathbb E}(X^{(j)})) \varepsilon}(u) du < \infty \quad \text{and} \quad 
  \sum_{k\geq 0} \int_0^{\alpha(k)}  Q^2_{(Y-{\mathbb E}(Y)) {e^{(j)}}}(u) du< \infty \, .
\end{equation}
For instance, \eqref{condalpha} is true  if, for some $q>2$,  $\sum_{k\geq 0} k^{2/(q-2)} \alpha(k) < \infty$ and
$${\mathbb E}(|(X^{(j)}-{\mathbb E}(X^{(j)})) \varepsilon|^q)< \infty \, , \quad
    {\mathbb E}(|(Y-{\mathbb E}(Y)) {e^{(j)}}|^q)< \infty \mbox{  for all }j \in \{1, \ldots , p\}\, . $$
    Note also that, if the sequence $(Y_i, X_{i}^{(1)}, \ldots ,X_{i}^{(p)})_{i \in {\mathbb Z}}$ is $m$-dependent (implying  that $\alpha(k)=0$ for $k>m$), then \eqref{condalpha} is exactly Assumption \ref{A_moment_3}.

Let us now described  how the matrix $\A$ of Proposition \ref{limitlaw} has to be changed in this dependent setting:
Let then $\A$ be the $2p\times 2p$ symmetric matrix defined as follows: 
\begin{itemize}
\item If $(j,k) \in \{1, \ldots, p\}^2$ then 
$$\A_{j,k}=\sum_{m \in {\mathbb Z}}{\mathbb E}((X_0^{(j)}-{\mathbb E}(X_0^{(j)}))(X_m^{(k)}-{\mathbb E}(X_m^{(k)}))\varepsilon_0 \varepsilon_m );$$ 
\item If $(j,k) \in \{p+1, \ldots, 2p\}^2$ then $$\A_{j,k}=\sum_{m \in {\mathbb Z}}{\mathbb E}((Y_0-{\mathbb E}(Y_0))(Y_m-{\mathbb E}(Y_m)) e^{(j-p)}_0 e^{(k-p)}_m);$$ \item  $(j,k) \in \{1, \ldots , p\} \times \{p+1, \ldots, 2p\}$ then $$\A_{j,k}=\sum_{m \in {\mathbb Z}}{\mathbb E}((X_0^{(j)}-{\mathbb E}(X_0^{(j)}))(Y_m-{\mathbb E}(Y_m)) \varepsilon_0 e^{(k-p)}_m).$$
\end{itemize}
One can easily check that the matrix $\A$ is well defined as soon as \eqref{condalpha} holds (more precisely each series described above converges absolutely under \eqref{condalpha}). 

Let $\B$ be the matrix defined before Proposition \ref{limitlaw}.
We are now in position to state the analogue of Proposition  \ref{limitlaw} in our $\alpha$-mixing context.

\begin{Proposition}\label{limitlawalpha}
Under Assumptions \ref{A_independent} and  \eqref{condalpha}
$$
  \sqrt{n}\Big([ \widehat{\bm{\alpha}}:\widehat{\bm{\theta}}]-
  [ {\bm{\alpha}}:{\bm{\theta}}]\Big)
  \quad \cvl
  {\mathcal N}_{2p}(0,{\mathbf B} {\mathbf A} {\mathbf B}).
$$
\end{Proposition}
\noindent{\bf Proof of Proposition \ref{limitlawalpha}.} It is the same as that of Proposition \ref{limitlaw}, by replacing the multivariate central limit theorem for i.i.d random vectors by the multivariate central limit theorem for $\alpha$-mixing random vectors (see Dedecker and Merlev\`ede \cite{DM2003}, Corollary 2 Item ($\alpha$)). \qed 

\medskip

As mentioned at the beginning of this section, the result of Proposition \ref{limitlawR2} (giving the asymptotic normality of $\sqrt n  ( \widehat R^2 -R^2 )$) remains true under the assumptions of Proposition \ref{limitlawalpha} (with the matrix $\A$ described above). The matrix $\B$ can be estimated as in Section \ref{Sec3} (using the ergodic theorem instead of the strong law of large numbers). Hence, to get a confidence interval for $R^2$, it remains to estimate the matrix $\A$. This is not an easy problem, and there is a large literature on this question. A possible way in the $\alpha$-mixing framework, is to use the general results on HAC (Heteroskedasticity and Autocorrelation Consistent) estimators stated in Theorem 1 of Andrews \cite{And91}. 

\medskip

To conclude, it should be noticed that the results of this section apply to the case of auto-regression, that is when 
$$
(Y_i, X_{i}^{(1)}, \ldots ,X_{i}^{(p)})=(X_i, X_{i-1},\ldots, X_{i-p})
$$
where $(X_i)_{i \in {\mathbb Z}}$ is a stationary sequence of $\alpha$-mixing random variables. To be more precise, let ${\mathcal M}_0= \sigma(X_i, i\leq 0)$ and ${\mathcal M}_k= \sigma (X_k, X_{k-1},\ldots, X_{k-p})$. Then  the coefficient $\alpha(k)$ defined in \eqref{defalpha} 
is exactly equal to the coefficient $\alpha({\mathcal M}_0, {\mathcal M}_k)$.

\section{Related results and applications} \label{Sec6}

In this Section, we give some additional results and present some possible applications of the techniques used in Section 3. In Section \ref{Sec4bis} we establish the asymptotic normality of the vector of individual $\widehat R^2$'s (from which one can  obtain confidence ellipsoids for the individual $R^2$'s).  In Section \ref{SecSA} we present an alternative to the usual sensitivity analysis.  In Section \ref{SecRS} we give a detailed outline for a robust screening method, which we illustrate with simulations when the output is a binary variable.

\subsection{Asymptotic joint distribution of individual $ \widehat R^2$'s}\label{Sec4bis}
\setcounter{equation}{0}
For $k \in \{1, \ldots , p\}$, let $U^{(k)}$ be the ${\mathbb L}^2$-subspace generated by ${\bf 1}$ and $X^{(k)}$, and let (assuming that $\text{Var}(Y)>0$)
$$
R^2_{(k)}=\frac{\text{Var}( \text{Proj}_{U^{(k)}}(Y))}{\text{Var}(Y)}\, .
$$
From  Lemma \ref{scalarproduct}, we have $R^2_{(k)}=\tau_k \theta_k$, where $\theta_k$ is defined by \eqref{def:theta} and 
\begin{equation}
    \label{def:tauk}
\tau_k=\frac{\mathrm{Cov}(Y, X^{(k)})}{\mathrm{Var}(X^{(k)})}.\end{equation} The empirical estimator of $R^2_{(k)}$ is then 
\begin{equation}
\label{def:tauhatk}
\widehat R^2_{(k)}= \widehat \tau_k \widehat \theta_k, \quad \text{where} \quad  \hat \tau_k = \frac{\sum_{i=1}^n (Y_i-\overline \Y)(X_{i}^{(k)}-\overline {\X}^{(k)})}{\sum_{i=1}^n (X_{i}^{(k)}-\overline {\X}^{(k)})^2}\, .
\end{equation}
The asymptotic distribution of $\sqrt n (\widehat R^2_{(k)}- R^2_{(k)})$ is given in Proposition \ref{limitlawR2}. In this section, we give the asymptotic joint distribution of 
$$
\sqrt n \left ((\widehat R^2_{(1)}\, \ldots, \widehat R^2_{(p)})^t- (R^2_{(1)}, \ldots, R^2_{(p)})^t \right ) \, .
$$

Let $\bm \tau = (\tau_1, \ldots, \tau_p)^t $ and $\widehat {\bm \tau} = (\widehat \tau_1, \ldots, \widehat \tau_p)^t$. The first step is to identify the limit distribution of 
$$
\sqrt n \left ( [\hat {\bm \tau} : \hat {\bm \theta} ]- [ \bm \tau :  \bm \theta ]\right ) \, .
$$
Let us then define the two matrices involved in this asymptotic distribution. For $k \in \{1, \ldots, p\}$, let 
$$\varepsilon^{(k)} = Y-\mbox{Proj}_{U^{(k)}}(Y) \quad \text{and} \quad \varepsilon^{(k)}_i = Y_i-\mbox{Proj}_{U^{(k)}_i}(Y_i),$$ where $U^{(k)}_i$  is the subspace of ${\mathbb L}^2$ generated by ${\bf 1}$ and $X^{(k)}_i$. Let also ${\bm \varepsilon}^{(k)}=(\varepsilon^{(k)}_1, \ldots, \varepsilon^{(k)}_n)^t$.
Recall that $e^{(k)}$ and $e^{(k)}_i$ have been defined in \eqref{def:ek}.

\begin{Assumption}
    \label{A_moment_4}
    ${\mathbb E}((X^{(j)}-{\mathbb E}(X^{(j)}))^2 {\varepsilon^{(j)}}^2)< \infty$ and 
    ${\mathbb E}((Y-{\mathbb E}(Y))^2 {e^{(j)}}^2)< \infty$ for all $j \in \{1, \ldots , p\}$.
\end{Assumption}
Let then $\C$ be the $2p\times 2p$ symmetric matrix defined as follows: 
\begin{itemize}
\item If $(j,k) \in \{1, \ldots, p\}^2$ then 
\begin{equation}
    \label{Cjj}
\C_{j,k}={\mathbb E}((X^{(j)}-{\mathbb E}(X^{(j)}))(X^{(k)}-{\mathbb E}(X^{(k)}))\varepsilon^{(j)}\varepsilon^{(k)});\end{equation}
\item If $(j,k) \in \{p+1, \ldots, 2p\}^2$ then $$\C_{j,k}={\mathbb E}((Y-{\mathbb E}(Y))^2 e^{(j-p)} e^{(k-p)});$$ \item  $(j,k) \in \{1, \ldots , p\} \times \{p+1, \ldots, 2p\}$ then $$\C_{j,k}={\mathbb E}((X^{(j)}-{\mathbb E}(X^{(j)}))(Y-{\mathbb E}(Y)) \varepsilon^{(j)} e^{(k-p})).$$
\end{itemize}

Let $\D$ be the $2p\times 2p$ diagonal matrix defined as follows :
\begin{itemize}
    \item 
    If $j \in \{1, \ldots, p\}$ then $\D_{j,j}= (\text{Var}(X^{(j)}))^{-1}$.
    \item 
    If $j \in \{p+1, \ldots, 2p\}$ then $\D_{j,j}= (\text{Var}(Y))^{-1}$.
\end{itemize}

\begin{Proposition}\label{limitlawVect}
Under Assumption \ref{A_moment_4} and assuming that
 $\mathrm{Var}(X^{(j)})>0$  
 for all $j \in \{1, \ldots , p\}$, we have
$$
  \sqrt{n}\Big([ \widehat{\bm{\tau}}:\widehat{\bm{\theta}}]-
  [ {\bm{\tau}}:{\bm{\theta}}]\Big)
  \quad \cvl
  {\mathcal N}_{2p}(0,{\mathbf D} {\mathbf C} {\mathbf D}).
$$
\end{Proposition}

\noindent{\bf Proof of Proposition \ref{limitlawVect}.} 
Let $\widehat \D$ be the $2p\times 2p$ diagonal matrix defined as follows :
\begin{itemize}
    \item 
    If $j \in \{1, \ldots, p\}$ then $\widehat \D_{j,j}= (\text{Var}_n(\X^{(j)}))^{-1}$.
    \item 
    If $j \in \{p+1, \ldots, 2p\}$ then $\widehat \D_{j,j}= (\text{Var}_n(\Y))^{-1}$.
\end{itemize}

Starting from \eqref{def:theta} and \eqref{def:alpha}, and noting that, for $k \in \{1, \ldots, p \}$,
\begin{align*}
  Y_i - \overline{\mathbf{Y}}&= \tau_k(X_{i}^{(k)}-\overline{\X}^{(k)}) + (\varepsilon^{(k)}_i- \overline{\bm{\varepsilon}}^{(k)})\, , \\
  (X_{i}^{(k)}- \overline {\mathbf{X}}^{(k)})&=\theta_k (Y_i - \overline{\mathbf{Y}}) +(e_{i}^{(k)}- \overline {\bm e}^{(k)})\, ,
\end{align*}
we see that (recall that ${\mathbb X}_0$ and $\tilde {\mathbb X}_0$ have been defined in \eqref{defX0} and \eqref{defX0tilde})
\begin{equation}\label{step1bis2}
\sqrt{n}\left([ \widehat{\bm{\tau}}:\widehat{\bm{\theta}}]-
  [ {\bm{\tau}}:{\bm{\theta}}]\right)= \widehat \D \frac{1}{\sqrt n} ({\bm{\varepsilon}^{(1)}}^t
{\mathbb X}_0^{(1)}, \ldots, {\bm{\varepsilon}^{(p)}}^t
{\mathbb X}_0^{(p)}, (\Y-\overline \Y{\bf 1}_n)^t {\bm e}^{(1)}, \ldots, (\Y-\overline \Y{\bf 1}_n)^t {\bm e}^{(p)})^t \, .
\end{equation}
One can easily check that
\begin{multline}\label{step2bis}
\frac{1}{\sqrt n} \Big(({\bm{\varepsilon}^{(1)}}^t
{\mathbb X}_0^{(1)}, \ldots, {\bm{\varepsilon}^{(p)}}^t
{\mathbb X}_0^{(p)}, (\Y-\overline \Y{\bf 1}_n)^t {\bm e}^{(1)}, \ldots, (\Y-\overline \Y{\bf 1}_n)^t {\bm e}^{(p)})^t \\-
({\bm{\varepsilon}^{(1)}}^t
\tilde {\mathbb X}_0^{(1)}, \ldots, {\bm{\varepsilon}^{(p)}}^t
\tilde {\mathbb X}_0^{(p)}, (\Y-{\mathbb E}(Y){\bf 1}_n)^t {\bm e}^{(1)}, \ldots, (\Y-{\mathbb E}(Y){\bf 1}_n)^t {\bm e}^{(p)})^t\Big )\cvp 0.
\end{multline}
Now, by the multivariate central limit theorem,
\begin{equation}\label{step3bis}
\frac{1}{\sqrt n} 
({\bm{\varepsilon}^{(1)}}^t
\tilde {\mathbb X}_0^{(1)}, \ldots, {\bm{\varepsilon}^{(p)}}^t
\tilde {\mathbb X}_0^{(p)}, (\Y-{\mathbb E}(Y){\bf 1}_n)^t {\bm e}^{(1)}, \ldots, (\Y-{\mathbb E}(Y){\bf 1}_n)^t {\bm e}^{(p)})^t \cvl \mathcal N_{2p}(0, \C).
\end{equation}
Moreover, by the strong law of large numbers,
\begin{equation}\label{step4bis}
    \lim_{n\rightarrow \infty} \widehat {\mathbf D} = {\mathbf D} \quad \text{almost surely}
\end{equation}
The result follows from \eqref{step1bis2}, \eqref{step2bis}, \eqref{step3bis} and \eqref{step4bis}. \hfill $\square$

\medskip

Let $\textbf{H}$ be the matrix with $p$ rows and $2p$ column such that
\begin{itemize}
    \item 
    If $(i,j) \in \{1, \ldots, p\}^2$ then $\textbf{H}_{i,j}= \theta_i \delta_{i,i}$.
    \item 
    If $ i \in \{1, \ldots, p\}$ and $j \in \{p+1, \ldots, 2p\}$ then $\textbf{H}_{i,j}=\tau_i \delta_{i,p+i} $.
\end{itemize}
As a consequence of Proposition \ref{limitlawVect}, we have 
\begin{Proposition}\label{limitlawR2ind}
Under Assumption \ref{A_moment_4} and assuming that
 $\mathrm{Var}(X^{(j)})>0$  
 for all $j \in \{1, \ldots , p\}$, we have 
$$
  \sqrt n \left ((\widehat R^2_{(1)}\, \ldots, \widehat R^2_{(p)})^t- (R^2_{(1)}, \ldots, R^2_{(p)})^t \right ) \cvl {\mathcal N}(0, \mathbf{H}\mathbf{D}\mathbf{C}\mathbf{D}\mathbf{H}^t).
$$
\end{Proposition}
\begin{Remark}
    Other quantities of interest are 
    $$
S^2_{(k)}=\frac{\mathrm{Var}( \mathrm{Proj}_{U^{(k)}}(Y))}{\mathrm{Var}(\mathrm{Proj}_{V}(Y))}\, \quad \text{for $k \in \{1, \ldots , p\}$},
$$
which are the proportions of the variance  of $\mathrm{Proj}_{V}(Y)$ that are explained by the best linear predictors in ${\mathbb L}^2$ based on the variables 
${\bf 1}, X^{(k)}$. Note that
$$S^2_{(k)}=\frac{R^2_{(k)}}{R^2}=\frac{\tau_k\theta_k}{{\bm{\theta}^t \bm{\alpha}}}\, .$$
Then, the asymptotic distribution of $\sqrt n \left ((\widehat S^2_{(1)}\, \ldots, \widehat S^2_{(p)})^t- (S^2_{(1)}, \ldots, S^2_{(p)})^t \right )$ may be derived from the asymptotic distribution of $ \sqrt n \left ( [\hat {\bm \tau} : \hat {\bm \theta} : \hat {\bm \alpha}]- [ \bm \tau :  \bm \theta : {\bm \alpha}]\right )$ via the delta method.
\end{Remark}
\noindent{\bf Proof of Proposition \ref{limitlawR2ind}.} 
Let $\phi : {\mathbb R}^{2p} \rightarrow {\mathbb R}^p$ be such that $\phi(x_1, \ldots , x_p, y_1, \ldots, y_p)=(x_1 y_1, \ldots , x_p y_p)^t$. Since 
$$
\left ((\widehat R^2_{(1)}\, \ldots, \widehat R^2_{(p)})^t- (R^2_{(1)}, \ldots, R^2_{(p)})^t \right ) = \phi( [ \widehat{\bm{\tau}}:\widehat{\bm{\theta}}]) - \phi([ {\bm{\tau}}:{\bm{\theta}}]) \, ,
$$
it suffices to apply the delta method to the function $\phi$. The results follows from Proposition \ref{limitlawVect} and the fact that the differential $D\phi_{\bm \tau,\bm \theta}$ of $\phi$ at point $(\bm\tau,\bm\theta)$ is given by 
$$
D\phi_{\bm \tau,\bm \theta}(h)=\textbf{H} h \, .\quad \quad \quad \quad \quad \square
$$

\medskip

To conclude this section, let us a give a consistent estimator of the matrix $\mathbf{H}\mathbf{D}\mathbf{C}\mathbf{D}\mathbf{H}^t$. We shall simply replace each matrix by its empirical counterpart.

Let then $\widehat {\mathbf C}$ be the $2p \times 2p$ symmetric matrix defined as follows: 
\begin{itemize}
\item If $(i,j) \in \{1, \ldots, p\}^2$ then 
$$\widehat {\mathbf C}_{j,k}=\frac 1 n \sum_{i=1}^n (X_{i}^{(j)}-\overline{\mathbf{X}}^{(j)})(X_{i}^{(k)}-\overline{\mathbf{X}}^{(k)})\widehat \varepsilon^{(j)}_i \widehat \varepsilon^{(k)}_i,$$ 
where $\widehat \varepsilon^{(k)}_i=(Y_i- \overline {\mathbf{Y}}) - \widehat \tau_{k} (X_i^{(k)} -\overline{\mathbf{X}}^{(k)})$;
\item If $(j,k) \in \{p+1, \ldots, 2p\}^2$ then $$\widehat {\mathbf C}_{j,k}= \frac 1 n \sum_{i=1}^n (Y_i-\overline{\mathbf{Y}})^2 \widehat e^{(j-p)}_{i} \widehat e^{(k-p)}_{i};$$ 
\item If $(j,k) \in \{1, \ldots , p\} \times \{p+1, \ldots, 2p\}$ then $$\widehat {\mathbf C}_{j,k}= \frac 1 n \sum_{i=1}^n(X_{i}^{(j)}- \overline {\mathbf{X}}^{(j)})(Y_i-\overline{\mathbf{Y}}) \widehat \varepsilon^{(j)}_i \widehat e^{(k-p)}_{i}.$$
\end{itemize}

\begin{Proposition}\label{estim2}
Assume that $\mathrm{Var}(X^{(j)})>0$, that ${\mathbb E}(Y^4)< \infty$ and that  ${\mathbb E}({X^{(j)}}^4)< \infty$  for any $j \in \{1, \ldots , p\}$. 
Then 
$$
\lim_{\nti} \widehat{\mathbf{H}}\widehat{\mathbf{D}}\widehat{\mathbf{C}}\widehat{\mathbf{D}}\widehat{\mathbf{H}}^t =\mathbf{H}\mathbf{D}\mathbf{C}\mathbf{D}\mathbf{H}^t \quad \mbox{almost surely.} 
$$
\end{Proposition}

\noindent{\bf Proof of Proposition \ref{estim2}.} By the strong law of large numbers, $\widehat {\textbf H}$ and $\widehat {\textbf D}$ converge almost surely to $\textbf{H}$ and $\textbf{D}$ respectively. Hence, it remains to prove that $\widehat \C$ converges to $\C$ almost surely, which can be done as in the proof of Proposition \ref{estim}.

\subsection{An alternative to the usual sensitivity analysis}\label{SecSA}

In sensitivity analysis, one usually assumes that
\begin{equation}
Y = m(X^{(1)},\ldots,X^{(p)}) 
\label{RegMod.eq} 
\end{equation}
(or that $Y = m(X^{(1)},\ldots,X^{(p)})+ \varepsilon$ as in \cite{SHMLT2017}) where $(X^{(1)},\ldots,X^{(p)})$ is a $p$ random vector with a known distribution $=
P_{1}\otimes\cdots\otimes P_{p}$, each $P_i$ having a moment of order two.
The number of variables $p$ may be large, and the function $m$  is an unknown function from $\R^{p}$ to
$\R$, which  may present strong non-linearities and high order interaction effects
between its coordinates.  On the basis of a $n$-sample $(Y_i, X_i^{(1)},\ldots,X_i^{(p)}), $ for $i=1,
\ldots, n$, we want to measure the impact of the main effects or interactions, i.e. to determine
the influence of each variable or group of variables on the output variable $Y$. 
The usual sensitivity analysis is mainly based on Hoeffding's decomposition, strongly related to the independency between 
the inputs variables $X^{(1)},\ldots,X^{(p)}$. In that latter case, if $m$
is  square integrable, one may
consider the classical Hoeffding-Sobol decomposition (Sobol~\cite{Sobol1993}) that leads to write
$m$ according to  its ANOVA functional expansion:
\begin{equation}
m(X^{(1)},\ldots,X^{(p)}) =
m_0+\sum_{{ \begin{array}{c}i_1<i_2<\cdots <i_k, \\i_i,\ldots,i_k\in\{1,\ldots,p\}\end{array}}} m_{i_1,\ldots,i_k} (X^{(i_1)},\ldots,X^{(i_k)})
\label{hoeffding.eq2}
\end{equation}
where the functions $m_{i_1,\ldots,i_k}$ are centered and orthogonal
in $\bL^{2}$ involving conditional expectation that is
\begin{eqnarray*}
m_0&=&\mathbb{E}(m(X^{(1)},\ldots,X^{(p)})), \quad m_j(X^{(j)})=\mathbb{E}(m(X^{(1)},\ldots,X^{(p)})|X^{(j)})-m_0\\
m_{j,k}(X^{(j)},X^{(k)})&=&\mathbb{E}(
m(X^{(1)},\ldots,X^{(p)})|X^{(j)},X^{(k)})-m_j(X^{(j)})-m_k(X^{(k)})-m_0, \quad \ldots
\end{eqnarray*}
Under the assumption  of independence on the $X^{(j)}$'s this decomposition is unique, with
$$\; \mathbb{E}(m_{i_1,\ldots,i_k}(X^{(i_1)},\ldots,X^{(i_k)})) = 0
$$
and for all $(i_1,\ldots,i_k) \neq  (i_1^\prime,\ldots,i_\ell^\prime)$,
$$\mathbb{E}\left (m_{i_1,\ldots,i_k}(X^{(i_1)},\ldots,X^{(i_k)})
m_{i_1^\prime,\ldots,i_\ell^\prime}(X^{(i_1^\prime)},\ldots,X^{(i_\ell^\prime)})
\right )=0 \, .
$$
This leads to the following variance decomposition:
$$
\Var(Y)
=
\sum_{{ \begin{array}{c}i_1<i_2<\cdots <i_k, \\i_i,\ldots,i_k\in\{1,\ldots,p\}\end{array}}} \Var \left(m_{i_1,\ldots,i_k} (X^{(i_1)},\ldots,X^{(i_k)})\right).$$

The Sobol 
sensitivity indices introduced by Sobol~\cite{Sobol1993} are
defined for any group $(X^{(i_1)},\ldots,X^{(i_k)})$ by
\begin{equation*}
 S_{i_1,\ldots,i_k} = \frac{\Var\left( m_{i_1,\ldots,i_k}(X^{(i_1)},\ldots,X^{(i_k)}) \right) }{\Var \left(Y\right)}.
\end{equation*}
They quantify
the contribution of a subset of variables $(X^{(1)},\ldots,X^{(p)})$ to the output
$Y$.
Several approaches are available for estimating these
sensitivity indices, see for example Iooss and
Lema\^itre~\cite{ioossLemaitre15} for a recent review. 

An alternative approach, which is much less computationally expensive, is as follows. Recall that the $X^{(i)}$'s are in ${\mathbb L}^2$ and let $\tilde X^{(i)} = X^{(i)}-{\mathbb E}(X^{(i)})$ (this centering may be done, since the distribution of $(X^{(1)},\ldots,X^{(p)})$ is known). Since the variables $X^{(j)}$'s are independent,  the family
$$
{\mathcal F}=\left \{ \tilde X^{(i_1)} \tilde X^{(i_2)}\cdots \tilde X^{(i_k)}, i_1 < \cdots < i_k\right \}
$$
is an orthogonal family of ${\mathbb L}^2$. We can then compute the individual $R^2$'s for this family, that is 
\begin{equation*}
 R^2_{i_1,\ldots,i_k} = \frac{\Var\left( \mbox{Proj}_{V^{(i_1, \ldots, i_k)}}(Y) \right) }{\Var \left(Y\right)}.
\end{equation*}
where $V^{(i_1, \ldots, i_k)}$ is the subspace of ${\mathbb L}^2$ generated by the product $\tilde X^{(i_1)} \tilde X^{(i_2)}\cdots \tilde X^{(i_k)}$. 

The main difference between the two approaches is that the sum of all the indices $R_{i_1,\ldots,i_k}$ is not equal to 1, unless $Y$ belongs to the subspace  $V_{\mathcal F}$ of ${\mathbb L}^2$ generated by the family ${\mathcal F}$. More precisely, we have that 
$$
\sum_{{ \begin{array}{c}i_1<i_2<\cdots <i_k, \\i_i,\ldots,i_k\in\{1,\ldots,p\}\end{array}}} R^2_{i_1,\ldots,i_k}= \frac{\Var\left( \mbox{Proj}_{V_{\mathcal F}}(Y) \right) }{\Var \left(Y\right)}.
$$
The advantage of using these $R^2$'s is that they are very easy to estimate (based on $n$ independent observations $(Y_i,X^{(1)}_i,\ldots,X^{(p)}_i)_{1 \leq i \leq n}$). In addition, as we have seen, we can give confidence intervals for these quantities.

Of course, the Sobol indexes and the $R^2$'s have not the same interpretation, the former are based on orthogonal complements obtained from conditional expectations (which are projections onto the spaces  ${\mathbb L}^2(\sigma(X^{(i_1)},X^{(i_2)}, \ldots, X^{(i_k)}))$), the latter are based on orthogonal projections onto the spaces $V^{(i_1, \ldots, i_k)}$ generated by the variables $\tilde X^{(i_1)} \tilde X^{(i_2)}\cdots \tilde X^{(i_k)}$.   For example, for a single index $i_k \in \{1, \ldots, p \} $, since $V_{i_k} \subseteq {\mathbb L}^2(\sigma(X^{(i_k)}))$, we always have $R^2_{i_k}\leq S_{i_k}$.





\subsection{A first step toward robust screening} \label{SecRS}
In this section, we outline a strategy for robust screening based on the results of Section \ref{Sec4bis}.

\subsubsection{Quick context} The basic idea of screening is to find, among a set of covariates $X^{(1)}, \ldots, X^{(p)}$, the ones having an association with an outcome $Y$. One usual way which is often cited consists in computing the $p$ absolute correlations $\widehat R_{(i)}$, as defined in Section \ref{Sec4bis}.
There are then two ways of thresholding these coefficients: either by keeping the largest $N_n$, or by keeping all those that exceed a threshold $\gamma_n$. 

To assess the performance of a screening method, whether it is based on computing $p$ correlations as we mentioned or not, one has to ask two questions:
\begin{enumerate}
    \item Is the method able to retrieve all the covariates that are associated with the outcome $Y$? This question is known in the literature as the Sure Screening Property  (see Fan and Lv \cite{Fan2008}).
    \item Is the method specific enough, i.e. can we ensure that not too many covariates are wrongly selected? This question directly relates to the control of the False Positive Rate (FPR).
\end{enumerate}
The questions of the sure screening property and the control of the FPR in a model-free environment is a growing field of research (see for instance Fan and Lv \cite{Fan2008}, Zhao and Li \cite{CoxScreening2012}, Fan and Song \cite{ScreeningGLM2010}, Pan et al. \cite{Panetal2019}).
To answer both questions, the choice of the threshold $N_n$ or  $\gamma_n$ is crucial.
One paper which is often cited is the paper by Fan and Lv \cite{Fan2008}. Their context is the following: they suppose a true joint linear model:
\begin{equation}\label{ModelFan}
    Y = \sum_{j=1}^p X^{(j)}\beta_j^\star + \varepsilon
\end{equation}
where $\varepsilon$ is a  Gaussian error.  
 The authors consider
the screening set: 
\begin{equation*}
    \widehat{\mathcal{M}} = \Big\{ j  : \widehat R_j \ \ \text{is among the $N_n$ largest}  \Big\}, 
\end{equation*}
where $N_n$ may be chosen as $N_n=[n/\log n]$. 
Ensuring the Sure Screening Property in that specific case consists in  showing that the support of $\beta^\star$ is \textit{at least} included in $\widehat{\mathcal{M}}$ with a probability close to one. 

Under certain conditions on the law of the $X^{(k)}$'s and assuming that the number of non-zero coefficients  in \eqref{ModelFan} is smaller than $N_n=[n/\log n]$, Fan and Lv \cite{Fan2008} are able to prove that the Sure Screening Property holds (see their Theorem 1). 
However, the threshold chosen  does not limit the number of false positives, i.e., the covariates that are not in the support of $\beta^\star$ but are marginally selected. This threshold, although not having any theoretical support for controlling the FPR, is widely used in the screening literature.

The two questions 1. and 2. above are addressed for instance in Zhao and Li \cite{CoxScreening2012} in the context of the Cox model. In this paper the authors propose to choose
\begin{equation*}
    \widehat{\mathcal{M}}= \left \{j : 
    \frac{\sqrt n |\widehat{\alpha}_j|}{\sqrt {\widehat v_j}} \geq \gamma \right \},
\end{equation*}
where $\widehat{\alpha}_j$ corresponds to the estimator in the marginal (and possibly misspecified) model, and  $\widehat v_j$ is a consistent estimator of the asymptotic variance of $\sqrt{n}(\widehat \alpha_j - \alpha_j)$.
According to Zhao and Li \cite{CoxScreening2012},  denoting by  $\Phi$  the distribution function of a standard Gaussian and by $q$  the expected false positive rate, the idea is to choose
\begin{equation*}
\gamma = \Phi^{-1}\left (1-\frac{q}{2} \right ).
\end{equation*}

In the screening literature, it is common to consider that independence between $X^{(j)}$ and $Y$ in the joint model is equivalent to being marginally associated with $Y$ (see Condition 3 in Fan and Lv~\cite{Fan2008} and Assumption 8 in Zhao and Li \cite{CoxScreening2012}). 
However, if we are working in a model-free context as in Section~\ref{Sec4bis}, these assumptions are unnecessary and it seems possible to define a threshold $\gamma$ using the distribution function $\Phi$.

\subsubsection{Our idea}
First, we highlight the fact that even in a model-free context, the question of screening the covariates is still relevant as a lowering dimension step before trying to fit a joint association between the remaining covariates and the outcome. 

Recall from Section~\ref{Sec4bis} the expression of $\tau_j$ and $\hat \tau_j$ given in Equation~\eqref{def:tauk} and \eqref{def:tauhatk}
\begin{equation*}
    \tau_j = \frac{\text{Cov}(Y,X^{(j)})}{\text{Var}(X^{(j)})} \ \text{ and } \ \widehat{\tau}_j = \frac{\sum_{i=1}^n(Y_i - \overline{\gY})(X_i^{(j)}-\overline{\gX}^{(j)})}{\sum_{i=1}^n(X_i^{(j)}-\overline{\gX}^{(j)})^2}.
\end{equation*}
Then, 
Proposition~\ref{limitlawVect} gives 
\begin{equation*}
    \sqrt{n}(\widehat{\tau}_j - \tau_j) \cvl \mathcal{N}(0,v_{j})
\end{equation*}
where
\begin{equation*}
    v_{j} =  \frac{\mathbb{E}((X^{(j)}-\mathbb{E}(X^{(j)}))^2{\varepsilon^{(j)}}^2)}{ (\text{Var}(X^{(j)}))^2}.
\end{equation*}
Moreover, it has been shown in Section~\ref{Sec4bis} that if we let
\begin{equation*}
    \widehat{v_j} = \frac{1}{n}  \sum_{i=1}^n \frac{\left (X_i^{(j)}-\overline{\X}^{(j)} \right )^2 
    \left( \widehat \varepsilon^{(j)}_i \right)^2}{\left (\text{Var}_n(X^{(j)})\right)^2}
\end{equation*}
then, \mbox{ for all } $j=1,\ldots,p$,
\begin{equation*}
    \lim_{n\rightarrow \infty} \widehat{v}_j  =  v_j \text{ almost surely} .
\end{equation*}
Similarly to Zhao and Li \cite{CoxScreening2012}, a reasonable screening set rule should be
\begin{equation}\label{Mhat}
     \widehat{\mathcal{M}}= \left \{j :\frac{\sqrt{n}|\widehat{\tau}_j|}{\sqrt{\widehat{v}_j }}\geq \gamma \right \} \ \mbox{ with }  \gamma = \Phi^{-1}\left(1-\frac{q}{2} \right).
\end{equation}
The next step is to show that this rule does indeed allow us to obtain the sure screening property, as well as satisfactory control of the FPR. We believe that moment assumptions on the variables should be sufficient to show these two properties, and we plan to study these issues in a future work. For now, in the subsection below, we provide a short simulation study, which shows that the procedure works well on a rather complex example.

\subsubsection{Simulation study}

For this simulation study the sample size $n$ will be equal to $n=500$, $n=1000$, $n=1500$ or $n=2000$,  with $1000$ explanatory variables, 14 being  correlated to the output $Y$ and 986 being not correlated to the output.

To highlight the robustness of our procedure, we will consider the case where the outputs $Y$ is a binary random variable. More precisely, the variable $Y$ will be generated by a logistic regression model through the first 10 explanatory random variables (which we call $design$). The variables $X^{(11)}, X^{(12)}, X^{(13)}, X^{(14)}$ are correlated to the design and to $Y$ but are not used to generate $Y$.

\paragraph{Description of the design.}
The design is such that


- $X^{(1)},X^{(2)},X^{(3)},X^{(4)},X^{(5)}$ are dependent centered Gaussian variables with covariance matrix given by the following Toeplitz matrix (with $\rho=0.57$)

$$\boldsymbol\Sigma=\begin{pmatrix}
1 & \rho & \rho^2 & \rho^3 & \rho^4 & \\
\rho & 1 &\rho &\rho^2 & \rho^3 \\
\rho^2 &\rho &1 &\rho &\rho^2& \\
\rho^3 &\rho^2 &\rho^1& 1 &\rho &\\ 
\rho^4 &\rho^3 &\rho^2& \rho& 1  \\

\end{pmatrix}$$

- $X^{(6)}\sim \mathcal{B}(0.35)$ is independent of the other variables of the design.

- $X^{(7)}\sim \chi^2(2)$ is independent of the other variables of the design.

- $X^{(8)}=X^{(1)} \times Z^{(1)}$, where $Z^{(1)} \sim \mathcal{P}(2)$ is independent of the other variables of the design.

- $X^{(9)}=X^{(2)}\times Z^{(2)}$,  where $Z^{(2)} \sim {\mathcal N}(1,1)$ is independent of the other variables of the design.

- $X_{10}\sim \mathcal{E}(1/2)$.
is independent of the other variables of the design.



\paragraph{Output.}

The output $Y$ is generated as follows:
$$ Y \sim \mathcal{B}\left (P\left (X^{(1)}, \ldots , X^{(10)}\right )\right )\, ,$$
with
$$
P\left (X^{(1)}, \ldots , X^{(10)}\right )=
\frac{1}{1+\exp(-\beta_0 - \beta_1 X^{(1)}- \cdots - \beta_{10} X^{(10)})}\, ,
$$
where $\beta=(\beta_0, \beta_1, \ldots, \beta_{10})$ is such that 
$$
\beta = \left (-1, 4 e^{-1/10}, -4 e^{-2/10}, 4 e^{-3/10}, -4e^{-4/10},2,4,6,3,3,4\right ) \, .
$$

\paragraph{Other explanatory variables.} The variables $X^{(11)}, X^{(12)}, X^{(13)}, X^{(14)}$ are as follows:

- $X^{(11)}=X^{(1)}+  Z^{(3)}$

- $X^{(12)}=X^{(3)}+  Z^{(4)}$

- $X^{(13)}=X^{(4)}+  Z^{(5)}$

- $X^{(14)}=X^{(5)}+  Z^{(6)}$

\noindent where $(Z^{(3)},  Z^{(4)},  Z^{(5)},   Z^{(6)})$ are independent of $(Y, X^{(1)}, \ldots, X^{(10)})$, and 
$Z^{(3)} \sim \mathcal{N}(0,0.3^2)$, $Z^{(3)} \sim \mathcal{N}(0,0.2^2)$, $Z^{(3)} \sim \mathcal{N}(0,0.35^2)$, $Z^{(3)} \sim \mathcal{N}(0,0.55^2)$.

The variables $X^{(15)}, X^{(16)}, \ldots, X^{(1000)}$ are independent of $(Y, X^{(1)}, \ldots, X^{(14)})$. These variables form a stationary Gaussian sequence with mean 0 and variance 1, and covariance matrix similar to the matrix $\boldsymbol\Sigma $ above (but now it is a 986 $\times$ 986 square matrix), with $\rho=0.7$.









\paragraph{The results.}
We apply the screening rule \eqref{Mhat}, with different values of $q$ (recall that $q$ represents the theoretical False Positive Rate). In our context a True Positive is an index between 1 and 14 selected by the procedure, since the first 14 variables are correlated with $Y$. A False Positive is an index between 15 and 1000 selected by the procedure. 

The results based on $N=1000$ repetitions for $q=15\%$ and different sample sizes, $n=500, 1000, 1500$ and $2000$ are given in Table \ref{tab7} below. This table give the rates $R_i$ of selection of the index $i$  in  $\{1, \ldots , 14\}$, the True Positive Rate (TPR) and the False Positive Rate (FPR). 

\begin{table}
\center
\begin{tabular}{|c|c|c|c|c|}
\hline
\quad \quad $n$ & 500 & 1000 & 1500 & 2000\\
\hline 
$R_1$ & 1 & 1 & 1 & 1\\
\hline
$R_2$ & 1 & 1 & 1 & 1\\
\hline
$R_3$ & 0.998 & 1 & 1 & 1\\
\hline
$R_4$ & 0.483 & 0.741 & 0.864 & 0.940\\
\hline
$R_5$ & 0.776 & 0.935 & 0.992 & 0.997\\
\hline
$R_6$ & 0.673 & 0.865 & 0.949 & 0.978 \\
\hline
$R_7$ & 1 & 1 & 1 & 1\\
\hline
$R_8$ & 1 & 1 & 1 & 1\\
\hline
$R_9$ & 1 & 1 & 1 & 1 \\
\hline
$R_{10}$ & 1 & 1 & 1 & 1\\
\hline
$R_{11}$ & 1 & 1 & 1 & 1 \\
\hline
$R_{12}$ & 0.997 & 1 & 1 & 1 \\
\hline
$R_{13}$ & 0.478 & 0.698 & 0.828 & 0.923\\
\hline
$R_{14}$ & 0.690 & 0.890 & 0.972 & 0.988\\
\hline
TPR & 0.864 & 0.938 &0.972 & 0.988 \\
\hline
FPR & 0.154 & 0.152 & 0.151 & 0.152\\
\hline
mean $ |\widehat{\mathcal{M}} |$ & 164.408 & 163.108 & 162.875 & 163.464\\
\hline
\end{tabular}
\caption{Case $q=15\%$, $n=500, 1000, 1500, 2000$. Selections rates for the indexes $i =1, \ldots , 14$, True Positive Rate, False Positive Rate, and mean $ |\widehat{\mathcal{M}} |$.}
\label{tab7}
\end{table}

The same results for $q=20\%$ are presented in Table \ref{tab8}.

\begin{table}
\center
\begin{tabular}{|c|c|c|c|c|}
\hline
\quad \quad $n$ & 500 & 1000 & 1500 & 2000\\
\hline 
$R_1$ & 1 & 1 & 1 & 1\\
\hline
$R_2$ & 1 & 1 & 1 & 1\\
\hline
$R_3$ & 0.997 & 1 & 1 & 1\\
\hline
$R_4$ & 0.579 & 0.798 & 0.913 & 0.949\\
\hline
$R_5$ & 0.783 & 0.956 & 0.990 & 0.998\\
\hline
$R_6$ & 0.689 & 0.894 & 0.964 & 0.988 \\
\hline
$R_7$ & 1 & 1 & 1 & 1\\
\hline
$R_8$ & 1 & 1 & 1 & 1\\
\hline
$R_9$ & 1 & 1 & 1 & 1 \\
\hline
$R_{10}$ & 1 & 1 & 1 & 1\\
\hline
$R_{11}$ & 1 & 1 & 1 & 1 \\
\hline
$R_{12}$ & 0.998 & 1 & 1 & 1 \\
\hline
$R_{13}$ & 0.549 & 0.762 & 0.888 & 0.928\\
\hline
$R_{14}$ & 0.723 & 0.904 & 0.971 & 0.990\\
\hline
TPR & 0.880 & 0.951 &0.980 & 0.990 \\
\hline
FPR & 0.205 & 0.202 & 0.202 & 0.202\\
\hline
mean $ |\widehat{\mathcal{M}} |$ & 214.161 & 212.603 & 212.976 & 212.928\\
\hline
\end{tabular}
\caption{Case $q=20\%$, $n=500, 1000, 1500, 2000$. Selections rates for the indexes $i =1, \ldots , 14$, True Positive Rate, False Positive Rate, and mean $ |\widehat{\mathcal{M}} |$.}
\label{tab8}
\end{table}

\paragraph{Concluding comments.} The results are quite satisfactory, even when $q=15\%$ (that is when we impose a small fraction of false positive). In that case the TPR increases from $86.4 \%$ when $n=500$ to $98.8\%$ when $n=2000$. We note also that 7 variables among the 14 that are correlated with $Y$ are systematically detected by the procedure when $n=500$, and 9 among 14 when $n\geq 1000$. As expected, the FPR is close (but a little larger) to $15\%$ when $q=15\%$, and to $20\%$ when $q=20\%$. The two objectives of the procedure (TPR close to 1 and FPR close to $q$) are thus achieved.

The advantage of this robust procedure is that we do not assume the existence of a true parametric relation between the output $Y$ and some of the explanatory variables. This dimension reduction procedure can therefore be applied in a wide variety of situations.

\bigskip

\noindent {\bf Acknowledgements.} We would like to thank the referees for their suggestions, 
which have helped to improve the content and presentation of this article.

\end{document}